\documentclass[12pt]{article}

 \usepackage{amssymb}

 \usepackage{amsfonts}

 \usepackage{amscd}

\usepackage{fancyhdr}
\usepackage{calc}
\newcommand{\be}{\begin{enumerate}}

\newcommand{\ee}{\end{enumerate}}

\let\ds=\displaystyle

\def\Inf{\mathop{\hbox{\rm inf}}}

\def\N{{\mathbb N}} \def\Z{{\mathbb Z}}

 \def\Q{{\mathbb Q}}

\def\R{{\mathbb R}} \def\C{{\mathbb C}}

 \def\F{{\mathbb F}}

\def\P{{\mathbb P}}

\def\s{{\bf s}}

\def\n{{\bf n}}
\def\l{{\bf {l}}}
\def\m{{\bf m}}
\def\i{{\bf i}}

\def\z{{\bf z}}
\def\e{{\bf e}}
\def\x{{\bf x}}
\def\nug{{\bf \nu}}
\def\a{{\bf a}}
\def\c{{\bf c}}
\def\r{{\bf r}}

\def\b{{\bf b}}
\def\m{{\bf m}}

\def\ord{{\rm {ord}}}

\def\eps{{ \varepsilon }}

\catcode`\ç=13
\defç{\c{c}}
\catcode`\é=13
\defé{\'e}
\catcode`\à=13
\defà{\`a}
\catcode`\è=13
\defè{\`e}
\catcode`\â=13
\defâ{\^a}
\catcode`\ù=13
\defù{\`u}
\catcode`\ê=13
\defê{\^e}
\catcode`\î=13
\defî{\^\i}
\catcode`\ô=13
\defô{\^o}

\newcommand{\la}{\langle}

\newcommand{\ra}{\rangle}

\newtheorem{Lemme}{Lemma}

\newtheorem{Corollaire}{Corollary}

\newtheorem{Theoreme}{Theorem}

\newtheorem{Definition}{Definition}

\begin{document}
\title{ Meromorphic Continuation of 
  Multivariable Euler Products and Applications.
\author{ Gautami Bhowmik\footnote{Universit\'e de Lille 1, UFR de
  Math., Laboratoire de Math. Paul Painlev\'e U.M.R. CNRS 8524,
  59655 Villeneuve d'Ascq Cedex, France.
Email : bhowmik@math.univ-lille1.fr}\hskip 0.4cm 
\& \hskip 0.3cm
Driss Essouabri\footnote{Universit\'e de Caen,
 UFR des Sciences, Campus 2,
 Laboratoire de Math. Nicolas Oresme, CNRS UMR 6139,
 Bd. Mal Juin, B.P. 5186, 14032 Caen, France.
Email : essoua@math.unicaen.fr}\hskip 0.4cm
\& \hskip 0.3cm
Ben Lichtin\footnote{Rochester  N.Y., USA. 
Email : lichtin@math.rochester.edu}
}}

\maketitle \noindent {\small {{\bf {Abstract.}} 
This article extends classical
one variable results about Euler products defined 
by integral valued polynomial or analytic functions 
to several variables.
We show there exists a  meromorphic continuation up 
to a presumed natural boundary, and also give a 
criterion, a la Estermann-Dahlquist, for 
the existence of a meromorphic extension to $\C^n.$
Among  applications we deduce analytic properties 
of height zeta functions for toric varieties over
$\Q$ and group zeta functions.

\medskip
{\bf Mathematics Subject Classifications:} 11M41, 11N37, 14G05, 32D15. \par
{\bf Key words: several variables zeta functions, Euler product, analytic
continuation, Manin's conjecture, Rational points, zeta functions of groups.} 
\vskip .2 in

{\bf Introduction:} 
\vskip .2 in

There are two fundamental problems in the study of Dirichlet series that admit an Euler product expansion in a region
of absolute convergence. The first problem is to prove  the existence of a meromorphic continuation into a larger region.
Assuming this is possible, the second problem is to describe precisely   the boundary of the domain  for
this meromorphic function. For Dirichlet series in one variable, the first important results go back to Esterman
\cite{estermann} who proved that if 
$h(Y)=\sum_d F(d)Y^d\,,$ where $F(d)$ is a ``ganzwertige" polynomial and $F(0)=1$,
then $Z(s)=\prod _p h(p^{-s})$ is absolutely convergent for $\Re(s) >1$ and can be meromorphically continued 
to the half plane $\Re(s) >0$. Moreover, $Z(s)$
be continued to the whole complex plane if and only if $h(Y)$ 
is a cyclotomic polynomial. 
Subsequently, Dahlquist \cite{dahlquist} extended this result to $h$  any analytic function with isolated singularities
within the unit circle.

The purpose of this paper is to extend these  two basic properties to a general class of Dirichlet series that have an
absolutely convergent  Euler product expansion in some open domain of $\C^n,\ n \ge 2.$  Thus, the object of our study is  an
Euler product
$$Z(h;\s) = \prod_p h(p^{-s_1},\dots, p^{-s_n}, p)$$ when
$h(X) = 1 + \sum_k h_k(X_1,\dots, X_n) X_{n+1}^k$ is either a polynomial or analytic function with integral coefficients. An
essential role in our analysis is played by a polyhedron in $\R^n,$ determined by the exponents of monomials appearing in the
expression for  $h(X).$ This polyhedron plays an important role in the Singularities literature, so it is, perhaps,  not too 
surprising to see it appear here as well. 

We first show that there  is a meromorphic continuation  up to a presumed natural boundary, whose geometry is
that of a tube over the boundary of a convex set. Our second main result  applies to
the case in which
$h$ depends only upon $X_1,\dots, X_n.$ In this event, we prove a very precise result that is the multivariate extension of the
work of Estermann-Dahlquist. This shows
that the presumed natural boundary is {\it the } natural boundary (in the sense given to this expression   in the
statement of Theorem \ref{ester+} in \S 1.2), unless
$h$ is a ``cyclotomic" polynomial. These results are proved in Section 1.

There are  several subjects, such as group theory, algebraic geometry, number theory,  knot theory, 
quantum groups, and combinatorics, in which  multivariate zeta functions can arise. Some of these are discussed in the
survey article of \cite{zagier}. It would therefore be interesting to find applications of our results/methods to the
analysis of such zeta functions. We discuss two applications in Sections 2, 3.
 
The first application (see Section 2)   originates with Manin's conjecture for toric varieties over $\Q.$  This gives a
precise description of the density of rational points with
``exponential height" at  most $t$ on such a variety.  Solutions  to
this conjecture   have been given by several authors (\cite{batyrev}, \cite{bretechetorique}, \cite{salberger})  (also see
\cite{moroz3}). In particular,  the method of de  la Bret\`eche used the deep work of Salberger to meromorphically continue a
certain generalized height zeta function into  some neighborhood of {\it exactly one  point} on the boundary of its domain of
analyticity. This   function was  a multivariate Dirichlet series with Euler product in the domain of absolute convergence.  
His approach sufficed to deduce the density asymptotic of interest for the conjecture, and also gave a strictly smaller
order (in $t$) error term. On the other hand, it did not address two general questions. The first inquires about all the other
points  on the boundary of the domain of analyticity of the Dirichlet series, in particular, how can they be
characterized/detected in general, or even calculated in concrete examples.  The second  asks  for an approximation to  the
natural boundary   of the meromorphic extension of the Dirichlet series. 

Describing precisely the entire   boundary of the domain of analyticity for this  series is needed to derive
the asymptotic  of rational points on the toric variety within a large family of expanding boxes. In the statement of Manin's
conjecture, only one expanding box appears, that  with sides all of the same length. There is however, no a priori reason why
this expanding box  should be privileged over any other.  Finding  the natural boundary of its meromorphic continuation appears
to be an interesting analytic problem by itself, and   has not, to our knowledge,  received any prior attention in the
literature. It may even encode something nontrivial    about the toric variety. One reason for believing this is the
observation that  the  boundary determines an estimate for the natural boundary of  a family of height zeta functions in one
complex variable. As such, it offers certain constraints upon the behavior of the zeroes and poles of each   height  function
in this family. Presumably, knowing something about such points, zeroes especially,  ought to be interesting.

We  solve   both of these problems,  using the   methods developed  in Section 1. As a result, our point of view is
rather different from that in the works cited above. Our main results are given in   Theorems \ref{application2},
\ref{application2bis},
\ref{coefficient} in
\S 2.3. Additional   discussion that contrasts our method and results with earlier work can be found   in
\S 2.3 following the statement of Theorem
\ref{application2}.  The last  result in \S 2 is  
Theorem \ref{application3} in  \S 2.4. This  
addresses a general (and natural) problem in the multiplicative
theory of integers, and is a good illustration of our method.  For {\it any} $n \ge 3,$ we give   the explicit 
asymptotic for the number of
$n-$fold products of positive integers that equal the $n^{th}$ power of an integer. The earlier papers (\cite{batyrev}, 
\cite{bretecheasterisque}, \cite{fouvry}, \cite{morozcubique}) had found the asymptotic when $n = 3.$ However, nothing
comparable   for arbitrary $n > 3$ seems to have been reported before in the literature.
 
The second application originates in group theory.  Several authors have associated  a Dirichlet series to certain
algebraic or finitely generated (nilpotent) groups  in order to study the density of finitely generated subgroups of large 
index. The algebraic structure of the groups that have been studied in this way enable the series to be written as an
absolutely convergent Euler product in one variable, whose factor at the prime  $p$ is an explicitly given function
$h(p, p^{-s}).$ In a series of papers,  du Sautoy, Grunewald and
others 
( \cite{dSG00}, \cite{dSG02}, \cite{dusautoy}) have
described with some success the analytic properties of  such  Euler products. 

The evidence produced in  these papers leads one to believe that  when there is ``uniformity" of the Euler product, there
should always exist a meromorphic extension, but that determining the natural boundary   is rather difficult in
general. On the other hand, the property of uniformity will not be satisfied for many other groups.   In this case,  only 
results that are less ambitious in nature  should be expected. For example, one can hope to study the boundary of analyticity
of the group's Dirichlet series. It is well known that this series has a real pole on this boundary.  A fundamental problem
had been  to show that this leading pole is rational, and that the series is meromorphic in some halfplane  that contains the
pole. The main result of
\cite{dSG00} established these two properties for any finitely generated nilpotent group. 

Our first observation in Section 3 is that the two main properties proved in \cite{dSG00}
can be established in a more elementary fashion, using Theorem 1 (see   \S 1.1) and certain diophantine estimates proved
in [ibid.]. Our second observation  is that the group zeta functions studied in  \cite{dSG02}, \cite{dusautoy},  
can be meromorphically extended outside a halfplane of absolute convergence by using the  method 
in Section 1. This is simpler than that used   in \cite{dusautoy}.  We also show that the presumed
natural boundary agrees with the one given in [ibid.].   The third observation addresses an analogous  problem about the
density of   subgroups   inside finite abelian groups of large order \cite{bhra1}.  We indicate by a simple example how 
nontrivial refinements of  standard density results can be found by using multivariate zeta functions and
Tauberian theorems.

A third example illustrates another way in which  the methods of this paper might eventually prove useful, but which we will 
not address further here. In the study of strings over
$p-$adic fields, one encounters Euler products in several variables. For example, in \cite{BFOW}, products of 5-point 
amplitudes  for the  ``open" strings are considered, where the amplitudes are defined as $p-$adic integrals
$$A_5^p(k_1,\dots, k_4)=\int_{\Q_p^2}\mid x\mid^{k_1k_2}\mid y\mid^{k_1k_3}\mid 1-x\mid^{k_2k_4}
\mid 1-y\mid^{k_3k_4}\mid x-y\mid^{k_2k_3} dxdy .  $$ 
The product $\prod_p A_5^p$ can be analytically continued. Indeed, our methods can certainly be used to prove this. In so
doing, one finds  interesting relations to the corresponding real amplitudes.

{\bf Notations:} \ For the reader's convenience, notations that will be used 
throughout the article are assembled here.
\be
\item 
  $\N=\{1,2,\dots \}$ denotes the set of positive integers,
$\N_0=\N \cup \{0\}$ and $p$ always denotes a  prime.
\item
The expression  \quad  $f(\lambda,{\bf y},{\bf x})\quad {\ll}_{{}_{{\bf y}}}\quad g({\bf x})\quad
 \mbox{ uniformly in }{\bf x}\in X \mbox{ and }{\lambda}\in \Lambda$ 
\newline means there exists $A=A({\bf y})>0$,
which depends neither on ${\bf x}$ nor  $\lambda$, but could eventually
 depend on  the   parameter vector  ${\bf y}$, such that:  $$
\qquad \forall {\bf x}\in X ~and~\forall {\lambda}\in
{\Lambda}\quad |f(\lambda,{\bf y},{\bf x})|\leq Ag({\bf x}). $$
When there is no ambiguity we  omit the word `uniformly'
above.
 \item For every ${\bf x}=(x_1,..,x_n) \in \R^n$, we set  
$\|{\bf x}\|=\sqrt{x_1^2+..+x_n^2}$ resp.  $|{\bf x}|=|x_1|+..+|x_n|$ to denote
the length resp. weight of $\bf x.$ We
denote the canonical basis of $\R^n$ by $(\e_1,\dots,\e_n).$ 
For every $\alpha=(\alpha_1,..,\alpha_n)\in \N_0^n\,,$ we also  set 
$\alpha!=\alpha_1!..\alpha_n!$. The standard inner product on $\R^n$ is denoted
$\la, \ra.$  
\item For every $s \in \C$, and for every non negative
$k$, we define $(_k^s)=\frac{s(s-1)..(s-k+1)}{k!}\,.$  
For two complex numbers $w$ and $z$, we define $w^z=e^{z\log w}\,,$
using the   principal branch  of the logarithm. 
We denote a vector in $\C^n$ by $\s=(s_1,\dots,s_n),$ and write
 $\s={\bold \sigma}+i{\bold \tau},$
 where ${\bold \sigma } = (\sigma_1,\dots,\sigma_n)$ and 
 ${\bold \tau }=(\tau_1,\dots,\tau_n)$ are the real resp. imaginary components of $\bold s$ (i.e. $\sigma_i=\Re(s_i)$ and
$\tau_i=\Im(s_i)$
  for each $i$). We also write $\la \x, \s\ra$ for $\sum_i x_i s_i$   if $\x \in \R^n, \s \in \C^n.$
\item The {\it unit polydisc}   $P(1)$ is 
the set $\{ \z=(z_1,\dots,z_n)\in \C^n\, :\  \sup_{i=1,\dots,n}
|z_i| < 1\}\,.$
\item Given $\alpha \in \N_0^n,$ we write $X^\alpha$ for the monomial $X_1^{\alpha_1} \cdots X_n^{\alpha_n}.$  For 
$h(X_1,\dots,X_n)= \sum_{\alpha \in \N_0^n} 
a_{\alpha} X^\alpha\,,$ the set $S(h):=\{\alpha \, :\  a_{\alpha} \neq 0\}$ is called the 
{\it support of} $h$. We also set $S^*(h):=S(h)\setminus \{0\}$. We denote by ${\cal E}(h)$ the boundary of the convex
hull of
$\bigcup \{\alpha +\R^n \, :\ \alpha \in 
S^*(h) \}$.  
This polyhedron is called the
 {\it Newton polyhedron} of  $h$ . We denote by
  $Ext (h)$ the set of extremal points of ${\cal
    E}(h)$ (a point of ${\cal E}(h)$ is  
extremal if it does not belong to the interior of any closed segment 
of ${\cal E}(h)$).
Obviously  $Ext (h)$ is a finite subset of $\N_0^n \setminus \{0\}$.
   
Similarly, if  $A \subset \N_0^{n}\setminus \{0\}$, 
we denote by ${\cal E}(A)$ the boundary of the convex hull of
$\bigcup \{\nug +\R_+^n \mid \nug \in A\}$ and call it the {\it Newton
polyhedron} of $A.$ Its set of extremal
points is denoted by $Ext(A)$.

\item If $A$ is a subset of $\N_0^n\setminus \{0\}$, we define 
${\widetilde {A}}$ 
as follows: 
\be
\item If $A$ is infinite set, then ${\widetilde {A}}$ denotes the 
set of $\nug \in A$ belonging to at least  one compact 
face of ${\cal E}(A).$ 
\item If $A$ is a finite set, then ${\widetilde {A}} = A.$
\ee
In either case, it is clear that ${\widetilde {A}}$ 
is a finite subset of $\N_0^n \setminus \{0\}.$  
The set ${\widetilde {A}}$ 
is called the saturation of $A.$
  
  \item Let ${\widetilde {A}}^o:
=\{\x \in \R_+^n \, :\  \forall \nug \in {\widetilde {A}},~~\la \x, \nug
\ra \geq 1\}$ be the dual of ${\widetilde {A}}$. Let $\iota (A)$ be the smallest weight of the elements of 
${\widetilde {A}}^o\,.$ We will call  $\iota (A)$ the index of $A$. We define 
$$R(A):= \{\alpha \in {\widetilde {A}}^o :
|\alpha| = \iota (A)\}.$$
For every $\alpha \in R(A)$, let  
$K(A;\alpha):=\{\nug \in {\widetilde {A}} \, :\  \la \alpha, \nug
\ra = 1\}$.
\ee

\section{Analytic properties of multivariate Euler products}
It will be convenient to split  the discussion in  two parts. The first   main result is Theorem 1. This constructs a
meromorphic extension for a large class of multivariate Euler products that converge absolutely in some product of
halfplanes of $\C^n.$ The second   main result, Theorem 2, extends the classical Estermann-Dahlquist   criterion for
the existence of a meromorphic extension to all of $\C^n, \ n \ge 2.$

\subsection{Meromorphic Continuation}
The first ingredient is the extension  of an Euler product, whose  $p^{th}$ factor $h(p^{-s_1},\dots, p^{-s_n})$
does not explicitly depend upon $p$ by itself, outside its domain of absolute convergence.  This extends Dahlquist's theorem
\cite{dahlquist} to several variables.   

\quad The following notations will be used. 
Let $\Lambda$ be an open subset of $\C^n$,  $l_1,\dots,l_r :\Lambda \to \C$ analytic functions,  and   $a_1,\dots,
a_r$   complex numbers. Define the Euler product 
$$Z_{\l}(\s)=Z_{\l}(s_1,\dots,s_n) =\prod_{p } \bigg(1+\sum_{k=1}^r
  \frac{a_k}{p^{l_k(\s)}}\bigg), $$
and for any $\delta \in \R$, set
$$W(\l;\delta)=W(l_1,\dots,l_r;\delta):=\{\s \in \Lambda : \forall i=1,\dots,r \quad \Re(l_i(\s))>\delta\}$$
It is clear that $s\mapsto Z_{\l}(\s)$ converges absolutely and defines a 
holomorphic function in the domain $W(\l;1)$.

\begin{Lemme}\label{zel}  

(i)  The function 
$Z_{\l}(\s)$  can be continued   into the domain $W(\l;0)$ as follows:  

there exists a set    $\{ \gamma (\n):  \n \in \N_0^r\} \subset \Q[a_0,\dots, a_r]$ 
  such that  for
every
$\delta >0,$    the function $G_\delta (\s)$ that is defined (and analytic) in $W(\l;1)$ by the equation 
$$ G_\delta (\s) = Z_{\l}(\s) \cdot \prod_{\n=(n_1,\dots,n_r)\in \N_0^r \atop
  1\leq |\n|  \leq [\delta^{-1}]}   
 \zeta \bigg( \sum_{j=1}^r ~ n_j l_j(\s) \bigg)^{-\gamma (\n)}$$
is actually a  bounded holomorphic function in $W(\l;\delta),$ 
 where it can be expressed  as an absolutely convergent Euler product.  

(ii)  When each $a_k \in \Z,$ each $\gamma (\n) \in \Z.$ In this case,  part (i) implies that the equation  
\begin{equation} \label{1} Z_{\l}(\s) = \prod_{\n=(n_1,\dots,n_r)\in \N_0^r \atop
  1\leq |\n|  \leq [\delta^{-1}]}   
 \zeta \bigg( \sum_{j=1}^r ~ n_j l_j(\s) \bigg)^{\gamma (\n)}\, G_\delta (\s) \end{equation}
determines a meromorphic extension of $Z_{\l}(\s)$ to $W(\l;\delta)$ for each $\delta > 0.$  
\end{Lemme} 

{\bf Remark.} As a result, only when each $\gamma (\n)$ is integral does it make sense to speak of a {\it meromorphic}
continuation of $Z_{\l}(\s)$ beyond $W(\l;1).$ For the
sake of simplicity, this function, defined by (1), in which each zeta factor  means, of course,  its meromorphic extension, 
is not   given a distinct notation. 

Even when this is not the case, part (i) shows that an analytic  extension of
$Z_{\l}(\s)$ is still possible in  simply connected subsets of any $W(\l;\delta),$ from which the branch (resp. polar) locus
of each    factor $\zeta (\n \cdot \l (\s))^{-\gamma (\n)}$ if $\gamma (\n) \notin \Z$ (resp. $\gamma (\n) \in \Z$) has been
deleted. For in each such subset, one can use the equation in (i) to  express $Z_{\l}(\s)$ as the product of $G_\delta (\s)$
with a single valued analytic continuation of each of the zeta factors. 
\hfill $\square$
 
{\bf Proof of Lemma \ref{zel}:} It suffices to prove part (i) since the proof of (ii) follows from the construction of the
$\gamma (\n)$ in (i).

Let $\delta \in (0,1)$ be arbitrary. To describe 
the   continuation of $Z_\l (\s)$  into  $W(\l;\delta),$ it will be convenient to work with a  somewhat larger class of Euler
products defined as follows: 

\begin{equation}\label{zlsg}
Z_\l (R_\delta;\s)=\prod_{p} \bigg(1+\sum_{k=1}^r
  \frac{a_k}{p^{l_k(\s)}}+ R_\delta(p;\s) \bigg) \quad
\end{equation}
where for all  $p$,
$s\mapsto R_\delta (p;\s)$ is a holomorphic function on $W(\l;\delta)$
satisfying\\ 
$ R_\delta (p;\s) \ll_{\l,\delta } p^{-2}$ uniformly in $p$
and   $\s \in W(\l;\delta)$.
Evidently,  $Z_\l (\s)=Z_\l (R_\delta;\s)$ when $R_\delta(p;\s)\equiv 0$. 

Now let us fix some notations: 
\be
\item     For each $m\in \N$, set
$${\cal L}_m(\l)
={\cal L}_m(l_1,\dots,l_r):=\{n_1 l_1+\dots+n_r l_r : n_1+\dots+n_r \geq m\};$$
\item   For each $\gamma_1,\dots,\gamma_r \in \C$, set 
\ $\Q_0[\gamma_1,\dots,\gamma_r] = \Z$ \ if \ $\gamma_1,\dots,\gamma_r \in \Z$
 and 
\newline 
\quad $\Q_0[\gamma_1,\dots,\gamma_r] = \Q[\gamma_1,\dots,\gamma_r]$ otherwise;
\item   $N=[2\delta^{-1}]\,;$
\item   $L(\s):= \prod_{k=1}^r {\zeta (l_k(\s))}^{-a_k}$ for $\s \in W(\l;1).$
\ee

By elementary computations, we obtain that for any  $\s \in W(\l;1):$ 
$$
L(\s) =
\prod_p \prod_{k=1}^r \bigg(1+\sum_{v_k=1}^N \frac{{a_k \choose v_k}
    (-1)^{v_k}}{p^{v_k l_k(\s)}}+H_N^k(p;\s)\bigg)$$
where, $\forall k=1,\dots,r$, \ $\s \mapsto H_N^k(p;\s)$ 
is a  holomorphic function  in $W(\l;\delta)$ and
satisfies the condition :
$H_N^k(p;\s)\ll_N p^{-\delta (N+1)} \ll_N p^{-2}$ uniformly
in $p$ and   $\s \in W(\l;\delta)$. It is also clear that   $a_k \in \N $ implies $H_N^k = 0$ once $N > a_k.$ 
\par
Thus, there exist 
$f_1,\dots,f_m \in {\cal L}_2(\l)$ and $d_1,\dots,d_m \in
\Q_0[a_1,\dots,a_r]$ such that :
$$
L(\s)= 
\prod_p \bigg(1-\sum_{k=1}^r \frac{a_k}{p^{l_k(\s)}}
+\sum_{i=1}^m \frac{d_i}{p^{f_i(\s)}}+K_N(p;\s)\bigg)
$$
 where $\s \mapsto K_N(p;\s)$ is a  holomorphic function  in $W(\l; \delta),$ satisfying 
the condition\\
 $K_N(p;\s)\ll_N p^{-2}$ uniformly
in $p$ and   $\s \in W(\l, \delta)$. \\
Now an easy computation shows that for every
 $\s \in W(\l, 1)$ :
\begin{eqnarray*}
Z_\l (R_\delta;\s) L(\s)&=& \prod_{p} \bigg(1+\sum_{k=1}^r
  \frac{a_k}{p^{l_k(\s)}} + R_\delta(p;\s)\bigg) \bigg(1-\sum_{k=1}^r \frac{a_k}{p^{l_k(\s)}}
+\sum_{i=1}^m \frac{d_i}{p^{f_i(\s)}}+K_N(p;\s)\bigg)\\
&=& \prod_{p} \bigg(1+\sum_{i=1}^m
  \frac{d_i}{p^{f_i(\s)}}-\sum_{k_1=1}^r\sum_{k_2=1}^r
  \frac{a_{k_1}a_{k_2}}{p^{l_{k_1}(\s)
+l_{k_2}(\s)}}
+\sum_{k=1}^r\sum_{i=1}^m \frac{a_k d_i}{p^{l_k(\s)+f_i(\s)}}+V_N(p;\s)\bigg)
\end{eqnarray*}
where $\s \mapsto V_N(p;\s)$ is a
  holomorphic function   in $W(\l;\delta),$ 
satisfying the bound\\ 
 $V_N(p;\s)\ll_N p^{-2}$  uniformly
in $p$ and  $\s \in W(\l; \delta)$. \par 
We have thus proved that there exist :
\be
\item 
$g_1,\dots,g_\mu \in {\cal L}_2(\l)$ and constants
$c_1,\dots,c_\mu \in \Q_0[a_1,\dots,a_r] $  
\item  for each   $p$   a
  holomorphic function \ $\s \mapsto R_{\delta,2}(p;\s)$    on
  $W(\l;\delta),$  satisfying \\ 
$R_{\delta ,2}(p;\s)\ll_\delta p^{-2}$ uniformly in $p$ and  $\s \in W(\l; \delta), $ 
\ee 
such that for every $\s \in
W(\l; 1)$ we have :
\begin{equation}\label{recu}
Z_\l (R_\delta;\s)\prod_{k=1}^r {\zeta (l_k(\s))}^{-a_k}=
\prod_{p} 
\bigg(1+\sum_{k=1}^\mu
  \frac{c_k}{p^{g_k(\s)}}+ R_{\delta ,2}(p;\s)\bigg).
\end{equation}
 
Since each  $g_k \in {\cal L}_2 (\l),$ it is clear that   $\Re\left(g_k(\s)\right) > 1$ 
for any $ \s \in W(\l; \frac{1}{2})$ and $  k=1,\dots,\mu.$  
This implies that for any $\delta'>max\left(\frac{1}{2},\delta\right)$:
 $$\s \mapsto \prod_{p} 
\bigg(1+\sum_{k=1}^\mu
  \frac{c_k}{p^{g_k(\s)}}+ R_{\delta ,2}(p;\s)\bigg)$$
is an absolutely convergent Euler product that is   holomorphic    in the domain $W(\l;\delta')$.\\

It is now evident how to proceed by induction.   Let $M= [\log_2 (N+1)]+1
\in
\N$. Repeating the above process  $M$  times, we conclude  that there exist: 
\be
\item  functions $h_1,\dots,h_q\in {\cal L}_1(\l)$ and constants
  $\gamma_1,\dots,\gamma_q \in \Q_0[a_1,\dots,a_r]$ 
\item 
 functions $u_1,\dots,u_\nu \in {\cal L}_{2^M}(\l)$ and constants
$b_1,\dots,b_\nu \in \Q_0[a_1,\dots,a_r]$ 
\item  for each $p,$  a  holomorphic function $\s \mapsto R_{\delta,M}(p;\s)$  on $W(\l;\delta),$ 
 satisfying \\
$R_{\delta ,M}(p;\s)\ll_\delta p^{-2}$
  uniformly in $p$ and  $\s \in W(\l; \delta)$ 
\ee 
such that for every $\s \in W(\l; 1)$ we have :
\begin{equation}\label{3}
 Z_\l (R_\delta; \s)\prod_{k=1}^q {\zeta (h_k(\s))}^{-\gamma_k}=\prod_{p} 
\bigg(1+\sum_{k=1}^\nu
  \frac{b_k}{p^{u_k(\s)}}+R_{\delta ,M}(p;\s) \bigg), 
\end{equation}
{\it and} the right side is absolutely convergent (and holomorphic) on $W(\l; \delta)$ since $2^{-M} < \delta/2.$ We now
multiply both sides of (\ref{3}) by $\prod_{\{k: h_k \in {\cal L}_{N+1}(\l)\}}\, \zeta (h_k(\s))^{\gamma_k}$ and  set 
$$G_\delta (\s): = Z_\l (R_\delta; \s) \cdot \bigg(\prod_{h_k \notin {\cal L}_{N+1}} {\zeta (h_k(\s))}^{-\gamma_k}\bigg).$$
In   $W(\l;1),$ \ $G_\delta (\s) = \prod_{\{k: h_k \in {\cal L}_{N+1}(\l)\}}\, \zeta (h_k(\s))^{\gamma_k} \cdot \prod_{p} 
\big(1+\sum_{k=1}^\nu
  \frac{b_k}{p^{u_k(\s)}}+R_{\delta ,M}(p;\s) \big).$ 
The preceding shows that the Euler product on the right is absolutely convergent in $W(\l;\delta).$ In addition, since $h_k\in
{\cal L}_{N+1}(\l)$ implies   $\Re\left(h_k(\s)\right)>(N+1)\delta  > 2,$ the product over $k$  also admits an analytic
continuation into $W(\l;\delta)$ as an absolutely convergent Euler product. Thus, $G_\delta (\s),$ whose individual factors
 in its definition   are, in general,   multivalued outside $W(\l;1)$ (with branch locus   the zero or polar
divisor of the individual factor), admits an analytic continuation into $W(\l;\delta)$ as an absolutely convergent Euler
product.   This proves (i). 

Part (ii) follows immediately from the fact that each $\gamma (\n)$ is
integral when each
$a_k$ is integral. Thus, the equation (1) determines a  meromorphic continuation of $Z_{\l}(\s)$ into $W(\l;\delta).$   

This  completes the proof of
Lemma \ref{zel}. 
\hfill $\square$  \par

\vskip .15 in

Let $h_0,\cdots,h_d$ be analytic functions on 
the unit polydisc $P(1)$ in $\C^n,$ satisfying  the  property $h_k(\bold 0)=0$ for each $k.$ 
Convergence in $P(1)$  should be understood as a normalization condition that can be easily
weakened without significant changes to the following discussion.  Define  
$$h(X_1,\dots,X_n,X_{n+1}) = 1 + \sum_{k=0}^d h_k(X_1,\dots,X_n)
X_{n+1}^k\,, \quad Z(h;\s)=\prod_{p} h(p^{-s_1},\dots,p^{-s_n},p)\,.$$ 
 Given the power series expansion of each $h_k,$\ \   
$h_k (X_1,\dots, X_n) = \sum_{\alpha \in \N_0^n} a_{\alpha, k} X^\alpha,$  

\indent  {\it   we assume throughout the rest of Section 1  that each $a_{\alpha, k} \in \Z.$} 

This will suffice for the applications of
interest in subsequent sections. 

\quad To state the first main result, we will also need the following notations, given the functions $h, h_0,\dots, h_d.$
\vskip .1 in

For each $\delta \in \R$, we set: 

$\qquad V(h; \delta):= \bigcap_{k=0}^{d}\big\{
\s \in \C^n :  \la \alpha, \bold \sigma \ra  >k+\delta \ \  \forall
\alpha \in Ext(h_k), {\mbox { and }}\   \sigma_i > \delta \ \ \forall i \big\}\,,$ 

and for $\delta > 0$ we set:   
\be
\item $N= \left[\frac{2(d+2)}{\delta}\right]+1$;
\item $ {\cal Y}_N: = \{(\alpha ,k) \in \N_0^n \times [0,
d] : 
\alpha \in S (h_k) {\mbox { and }} 1\leq |\alpha |\leq N\}\,, \ \, r_N:= \# {\cal Y}_N,$ and 
${\cal N}(\delta):=\{\n = (n_{\alpha,k}) \in \N_0^{r_N} :
1\leq |\n|  \leq \delta^{-1}\}\,.$ 
\ee

\begin{Theoreme}\label{zeh} 
There exists $A>0$ such that $Z(h;\s)$ converges absolutely in
$V(h;A)$. In addition, $Z(h;\s)$ can be continued into the
domain $V(h;0)$ as a meromorphic function as follows. 
For any $\delta > 0,$ there exists   $\{ \gamma (\n):
\n
\in {\cal N}(\delta)\}
\subset
\Z$ and $G_\delta (\s),$ a bounded holomorphic function on $V(h;\delta),$  such that the equation  
\begin{equation} \label {(4)} Z(h;\s)  = \prod_{\n=(n_{\alpha,k}) \in {\cal N}(\delta)}
\ 
{ \zeta \bigg(\sum_{(\alpha ,k) \in {\cal Y}_N}  n_{\alpha,k}
\left(\la \alpha , \s\ra -k\right)\bigg)}^{\gamma (\n)}~ \cdot G_\delta (\s),
\end{equation}
a priori valid  in $V(h;A),$ extends to $V(h;\delta)$ outside the polar divisor of the product
over $\n \in {\cal N}(\delta).$
Moreover $G_\delta$ can be expressed as an absolutely convergent Euler
product in the domain $V(h;\delta)$.
\end{Theoreme}

{\bf Proof:} The idea is to reduce the problem  to that studied in Lemma 1. The first
needed observation is clear. 
\begin{Lemme}
Let $k\in\{0,\dots,d\},$ and $\s  \in \C^n$ be such that 
 $\bold \sigma \in (0, \infty)^n$. Then 
\newline 
$\ds \Inf_{\alpha \in  S(h_k)}  \la \alpha , \bold \sigma \ra  = 
\Inf_{\alpha \in Ext (h_k)}   \la \alpha , \bold \sigma \ra $.
\end{Lemme}

\medskip
Next, we fix $A=2 (1+d) +1,$ and let $\s \in V(h;A).$
Evidently, this implies $\la \alpha, \bold \sigma \ra  \geq A |\alpha |$ for all $\alpha \neq \bold 0,\  \alpha \in
\cup_{k=1}^d  S(h_k).$   Thus for each $k\in [0, d]$, the convergence of $h_k$ on   $P(1)$   implies: 
 \begin{eqnarray*}
\sum_{\alpha \in S(h_k)}\left|\frac{a_{\alpha,k}}{p^{\la \alpha, \s \ra -k}}\right| 
& \leq & 
\sum_{\alpha \in S(h_k)}\frac{|a_{\alpha, k}|}{p^{  \la \alpha, \bold \sigma \ra - k}}
\leq 
\sum_{\alpha \in S(h_k)}
\frac{|a_{\alpha,k}|}{p^{A |\alpha |- k}} \\
&\leq &
\sum_{\alpha \in S(h_k)}
\frac{|a_{\alpha,k}|}{2^{ A |\alpha |/2 }} \cdot \frac{1}{p^{\frac{A}{2}- k}}
\ll \frac{1}{p^{\frac{A}{2}- k}} \ll \frac{1}{p^{\frac{A}{2}- d}}.
\end{eqnarray*}
By the definition of $A,$ we have $ \frac{A}{2}- d >1$. 
We conclude that
$$\s \mapsto Z(h;\s)=\prod_{p} h(p^{-s_1},\dots,p^{-s_n},p)
= \prod_{p } \bigg(1+\sum_{k=0}^d \ \sum_{\alpha \in S(h_k)}\frac{a_{\alpha,k}}{p^{\la \alpha, \s \ra -k}}\bigg)$$
is a  holomorphic function
in   $V(h;A)$. For any $\delta \in (0, A)$ and $\s \in V(h;\delta),$ it is then easy to express $Z(h;\s)$ in the form of (1) 
by subtracting off a sufficiently large tail of each $h_k$ that depends upon $\delta.$ The details are as follows.

As above,  $N (= N_\delta) = [2(d+2) \delta^{-1}]+1.$   For   $\s \in V(h;\delta)$ note that 
\begin{eqnarray*}
\sum_{k=0}^d \sum_{\alpha \in S(h_k) \atop |\alpha | \geq N}
\left|\frac{a_{\alpha,k}}{p^{\la \alpha, \s \ra -k}}\right| 
\leq 
\sum_{k=0}^d \sum_{\alpha \in S(h_k) \atop |\alpha | \geq N}
\frac{|a_{\alpha,k}|}{p^{ \la \alpha, \bold \sigma \ra  - k}} 
& \leq & 
\sum_{k=0}^d \sum_{\alpha \in S(h_k) \atop |\alpha | \geq N}
\frac{|a_{\alpha,k}|}{p^{\delta |\alpha |- k}} \\
&\leq &
\sum_{k=0}^d \frac{1}{p^{\frac{\delta 
      N}{2}- k}} \cdot \sum_{\alpha \in S(h_k)}
\frac{|a_{\alpha,k}|}{2^{\delta |\alpha| / 2}} 
\ll \sum_{k=0}^d \frac{1}{p^{\frac{\delta N}{2}- k}} \ll
\frac{1}{p^{\frac{\delta N}{2}- d}}.
\end{eqnarray*}
Since $\frac{\delta N}{2}- d > 2$, we conclude that $Z(h;\s)$ can be rewritten for any $\s \in V(h;A)$  as  
$$Z(h;s)=\prod_{p } \bigg(1+\sum_{k=0}^d \ \sum_{\alpha \in
    S(h_k) \atop 1\leq  |\alpha |\leq N}\frac{a_{\alpha,k}}{p^{\la \alpha,
      \s \ra -k}} + R_\delta (p;\s)\bigg),$$
where $\s \mapsto R_\delta (p;\s)$ is a bounded holomorphic function in
$V(h;\delta)$ such that 
\newline 
$R_\delta (p;\s) \ll_\delta p^{-2}$\  
uniformly in $p$ and  $\s \in V(h;\delta)$. \par 
The procedure described in Lemma 1 then applies with the
finite set of functions    $\s \to  \la \alpha, \s \ra - k,$ when
$1\le |\alpha|\le N, \ 0 \le k \le d.$ Thus,   for any
$\delta \in (0, A),$
$Z(h;\s)$ can   be analytically continued as a meromorphic function in  
$V(h;\delta),$  whose  precise expression is given   by  (5).
\hfill $\square$

{\bf Remark:}  The preceding argument actually shows that $Z(h;\s)$ converges absolutely in $V(h;1),$ even if $1 < A.$ 
This observation will be needed in the proof of Theorem 3 (see \S 2.1). The details justifying  this assertion are as follows.
For any   $\delta >0,$  the preceding discussion has shown that 
$Z(h;\s)\big|_{V(h;A)}$ can be rewritten   as
\begin{equation}\label{remeqn} Z(h;s)=\prod_{p } \bigg(1+\sum_{k=0}^d
\ \sum_{\alpha \in S (h_k) \atop
1\leq  |\alpha |\leq N}
\frac{a_{\alpha,k}}{p^{\la \alpha,\s \ra -k}} +
R_\delta (p;\s)\bigg),\end{equation}
where $N = N_\delta$ is defined as above, and $\s \mapsto R_\delta (p;\s)$ is a
bounded holomorphic function in
$V(h;\delta)$ that satisfies the bound 
$R_\delta (p;\s) \ll_\delta p^{-2}$
 uniformly in $p$ and  $\s \in V(h;\delta)$.

Thus, if $\delta = 1 < A,$ then 
  $\s \mapsto R_1 (p;\s)$ is a bounded
holomorphic function in
$V(h;1)$ such that
   $R_1 (p;\s) \ll_\delta p^{-2}$
while the sum of the finitely many terms indexed by those $\alpha \in S(h_k)$ and $|\alpha| \le N_1$ satisfies the property
that for any compact subset ${\cal K} \subset V(h;1),$ there exists $\theta_{{\cal K}} > 0$ such that the sum is 
$O(p^{-1-\theta_{{\cal K}}})$ uniformly in  ${\cal K}.$  Thus, the product in (\ref{remeqn})  converges absolutely in
$V(h;1).$\hfill
$\square$
 
\vskip .1 in

A simple extension of Theorem 1 will also be useful for the discussion in Section 3. This enlarges the original domain  from
which one begins the meromorphic extension of $Z(h;\s),$   by allowing   some 
$\sigma_i$ to be negative. The case when $h$ is a polynomial is the most naturally occuring one, so it will be given 
below. A simple extension to allow suitable rational factors in $X_{n+1}$ can also be made, but this need not be done here. 
As above, we write  $h = 1 + \sum_{k=0}^{d} \ h_k(X_1,\dots, X_n) X_{n+1}^k  \in
\Z[X_1,\dots, X_n, X_{n+1}],$ and assume   $h_k (\bold 0) = 0$ for each $k.$   Set:  
$$h_k = \sum_{\alpha \neq \bold 0} a_{\alpha, k} X_1^{\alpha_1} \dots X_n^{\alpha_n},  \ \ \l = (l_{\alpha, k})_{(\alpha,
k)}\,,\ \ {\mbox { where }}\ l_{\alpha, k}(\s) = \la \alpha, \s \ra - k  \ \ {\mbox { iff }}\  \alpha \in S(h_k).$$     
For any $ \delta \in \R$, set
$$\ds V^\# (h;\delta):= \bigcap_{k=0}^{d}\bigg\{
\s \in \C^n :  \la \alpha, \bold \sigma \ra  >k+\delta \ \  \forall
\alpha \in Ext(h_k)   \bigg\}\,.$$
The  proof of the following assertion is now straightforward.
\begin{Corollaire}\label{zehf}   $\s \mapsto Z(h;\s) $ 
can be continued meromorphically  from   $V^\# (h;1)$ (where $Z(h;\s)$   converges absolutely), into   $V^\#(h;\delta).$ 
\end{Corollaire}

{\bf Proof:} Apply the proof of Lemma \ref{zel} using the map $\l,$ as above. It is clear that for any $\delta,$ \ $\s \in 
W(\l;\delta)$ if and only if \ $\s \in V^\# (h; \delta).$ Thus,  the expression for the meromorphic continuation of
$Z(h;\s)$ in
$V^\#(h;\delta)$ follows directly from (1).
 
\subsection{ The natural boundary} 

In this subsection we work with a single analytic function \ $h(X) = 1 + \sum_{\alpha \neq \bold 0} a_\alpha 
X^\alpha.$ 
In the setting of Theorem 1, one thinks of $h$ as the function denoted $1 + h_0.$ Thus,  
 $$V(h;\delta): =\{
\s \in \C^n :  \la \alpha, \bold \sigma \ra  >\delta  \quad \forall
\alpha \in  Ext (h) \ {\mbox { and }}\   \sigma_i > \delta \ \ \forall i \big\}\,.$$  

The second main result of \S 1 concerns the Euler product   
$$ Z(h;\s) = 
\prod_{p } h(p^{-s_1},\dots,p^{-s_n}) =
\prod_{p} \bigg(1 + \sum_{\alpha \neq \bold 0}
  \frac{a_\alpha }{p^{\la \alpha,\s \ra}}\bigg).$$
Theorem 1 has shown that $Z(h;\s)$ can be meromorphically continued to $V(h;0)$ from some domain $V(h;A),\ A > 1, $ where
it converges absolutely as an Euler product. Of interest then are conditions satisfied by $h$ that imply  $Z(h;\s)$ can
or cannot be extended still further.

\begin{Theoreme}\label{ester+}
Assume each $a_\alpha \in \Z,$ and there exist $C,D>0$ such that for all $\alpha \in
\N_0^n$,
  $$|a_\alpha| \leq C (1+| \alpha |)^D \,.$$ Then
$Z(h;\s)$ can be continued to   $\C^n$ as a meromorphic function 
if and only if $h$ is `cyclotomic', i.e. there exists
a finite set $(\m_j)_{j=1}^q$ of elements of 
$\N_0^n \setminus \{\bold 0\}$ and a finite set of
integers $\{\gamma_j = -\gamma (\m_j)\}_{j=1}^q$ such that: 
$$h(X)=\prod_{j=1}^q (1-X^{\m_j})^{\gamma_j}= 
\prod_{j=1}^q (1-X_1^{m_{1,j}}\dots X_n^{m_{n,j}})^{\gamma_j}.$$
In all other cases the boundary   $\partial V(h;0)$ is the natural boundary. For purposes of this paper, this expression 
means that  
$Z(h;\s)$ can not be continued meromorphically into $V(h;\delta)$ for any $\delta < 0.$ 
\end{Theoreme}

\vskip .1 in

{\bf Proof:} It is clear that if $h$ is cyclotomic then $Z(h;\s)$ has a meromorphic extension to $\C^n.$ So, it suffices
to prove the converse. To do so, it suffices to assume only that $Z(h;\s)$ admits a meromorphic extension to $V(h;\delta_0)$
for some $\delta_0 < 0.$ The argument to follow will then show  that $h$ must  be cyclotomic, from which it follows immediately
that  $Z(h;\s)$ is   meromorphically extendible to $\C^n.$

We denote the elements of $Ext (h)$ by setting 
$Ext (h) = \{\alpha_1,\dots,\alpha_q\}\,.$ 
 
By the proof of Theorem 1, the continuation of $Z(h;\s)$ into  each  $V(h;\frac 1r), r = 1, 2,  \dots$ is determined by the
following property. There exist  $A \ge 1,$   a sequence $\{\gamma (\m)\}_{\m \in \N_0^n}$ of integers, and  a strictly
increasing sequence of positive integers $\{N_r\}_r$ such that  for each  
$\s \in V(h;A)$ and $r \ge 1:$
\begin{equation}\label{gdeltaa} 
  Z(h;\s) =\bigg(\prod_{m\in \N_0^n \atop 1\leq
   |\m|\leq N_r } {\zeta\left(\la \m, \s \ra \right)}^{\gamma(\m)}\bigg)
\times G_{1/r}(\s),
\end{equation}
where $G_{1/r}(\s)$ 
 is an absolutely convergent Euler product that is  bounded and
holomorphic in  $V(h;\frac{1}{r}).$ Thus, the extension of $Z(h;\s)$ into $V(h;\frac1r)$ is given explicitly as a product of 
$G_{1/r}$ with  the
meromorphic continuation into this domain of each of the finitely many zeta factors in (\ref {gdeltaa}).  

Set ${\it Ex} :=\{\m \in \N_0^n \setminus \{\bold 0\} : \gamma(\m) \neq 0\}$ \ and \ 
${\it  Ex }_-:=\{\m \in \N_0^n \setminus \{\bold 0\} : \gamma(\m) <0\}$. \par
We have to distinguish two cases:\par
{\bf {  {Case 1: ${\it  Ex }$ is infinite}}} \par
As above, let   $\delta_0 <0$ be such that $Z(h;\s)$ has a 
meromorphic continuation to $V(h;\delta_0)$.\\
Let $\rho_0$ be any fixed (and necessarily nonreal) zero of the Riemann zeta function satisfying $\Re(\rho_0)=\frac{1}{2}$.\\
Fix $\beta =(\beta_1,\dots,\beta_n) \in (0, \infty)^n$ such that 
$\beta_1,\dots,\beta_n$ are $\Q$-linearly independent, and set $Z_\beta(t):= Z(h; t\beta)$.\par
For all $\m \in {\it  Ex  }$ we set 
$t_\m = \frac{1}{\la \m, \beta \ra}$  if $\gamma(\m)<0,$  \  
and \ $t_\m = \frac{\rho_0}{\la \m, \beta \ra}$  if $\gamma(\m)>0$.\\ 
In addition, choose   for each $\m \in K$, $r(\m) \in \N$ satisfying: 
$$\ds r(\m) > \frac{ 2 \cdot |\m| \cdot \sup_i \beta_i}
{\inf_j\, \la \alpha_j , \beta \ra} \qquad \mbox{ and} \qquad r(\m)\geq |\m|.$$ 
It follows that $N_{r(\m)} \geq r(\m) \geq |\m|.$ By   (\ref{gdeltaa}), we have for each $\m \in {\it  Ex  }$ and $t\beta
\in V(h;
\frac 1{r(\m)}):$
\begin{equation}\label{gdeltam}
Z_\beta (t)=Z(h;t\beta)=\zeta(t \la \m, \beta \ra )^{\gamma(\m)}
\bigg(\prod_{\m'\in \N_0^n \setminus \{\m\} \atop 1\leq
   |\m'|\leq {N_{r(\m)}} } \zeta(t\la \m', \beta \ra )^{\gamma(\m')}\bigg)
 G_{1/r(\m)} (t\beta).
\end{equation}

From the definition of $r(\m)$, it follows that for each $\alpha_j \in Ext (h)$:
$$\Re(\la \alpha_j, t_\m \beta  \ra) \geq \frac{\la  \alpha_j, \beta\ra}
{2 \cdot \la \m, \beta \ra }\geq \frac{\la \alpha_j, \beta \ra}
{ 2 \cdot |\m| \cdot \sup_i \beta_i }>\frac{1}{r(\m)}.$$
Thus, $t\mapsto G_{1/r(\m)} (t\beta)$ is  holomorphic 
 in a neighbourhood of $t = t_\m $.\par
We now distinguish two subcases: \par

\vskip .1 in

{\bf {First subcase: ${\it  Ex }_-$ is infinite}} \par
Let $\m \in {\it  Ex }_-,$ so that   $t_\m = \frac{1}{\la \m, \beta \ra} > 0.$ It follows that $t_{\m}$ is not a pole of
$\zeta(t\la \m', \beta \ra )^{\gamma(\m')}$ for every $\m' \neq \m \in
\N_0^n\,.$    This is clear if $\gamma (\m') > 0$  since 
the only possible  pole of this function occurs when 
$t = \frac{1}{\la \m', \beta \ra},$ which cannot equal $t_{\m}$ because  
$t_{\m} =\frac{1}{\la \m, \beta \ra} \neq \frac{1}{\la \m', \beta \ra}$. If $\gamma (\m') < 0,$ then poles of $\zeta(t \la \m',
\beta \ra )^{\gamma(\m')}$ must be zeroes of $\zeta(t \la \m', \beta \ra ).$ A classical fact (\cite{tit}, pg. 30)
tells us that there are no positive zeroes of
$\zeta (s).$ Thus, $t_{\m}$ cannot  be a pole of $\zeta(t \la \m', \beta \ra)^{\gamma(\m')}.$   On the
other hand,  
$\gamma (\m) < 0$ implies that $t_\m$ \, {\it is}  a zero of $Z_\beta (t)$ since  $|\m| \le N_{r(\m)}.$
 
Furthermore, it is clear that the sequence $\{t_{\m}\}_{\m \in {\it  Ex}_-}$ of zeroes of $Z_\beta (t)$ converges to
$0$ when
$|\m|
\rightarrow +\infty.$

Now, if $Z(h;\s)$ had a meromorphic continuation to $V(h; \delta_0),$  then  $Z_\beta (t)$ would have to have a
meromorphic continuation to $U(\delta_1):=\{ t
\in \C : \Re(t) >\delta_1 \},$ where 
$\delta_1 = \sup_{1\leq j\leq q} \frac{\delta_0}{\la \alpha_j , \beta
  \ra} < 0$. Thus, $Z_\beta (t)$ would have to be identically zero, which is impossible because each $G_{1/r}(\s)$ is an
absolutely convergent Euler product in $V(h;1/r),$ and cannot therefore be identically zero.    We conclude that in this
subcase,  
$Z(h;\s)$ cannot be meromorphically extended to any $V(h;\delta)$ when $\delta < 0$.\par
\vskip .1 in

{\bf {Second  subcase:   ${\it  Ex }_-$ is finite}} \par
Choose  $a>0$ such that $\zeta(z) \neq 0$ for  $|z| \leq a$. \\
Set 
$$  B:= 2 \cdot \frac{(\sup_i \beta_i) \cdot |\rho_0| \cdot 
  (\sup_{\m\in {\it  Ex }_-} |\m|)}{a \cdot (\inf_i \beta_i)} >0 .$$
Define ${\it  Ex }_+: = {\it  Ex } \setminus {\it  Ex }_-\,,$ and fix $\m \in {\it Ex}_+$ such that   $|\m| \geq B$. 
Then $\gamma(\m) >0$ and $t_\m = \frac{\rho_0}{\la \m, \beta \ra}
\in \C \setminus \R$. \\
We then observe the following:
\be
\item
for all $\m' \in {\it  Ex }_+$ satisfying $\m' \neq \m$, $t_{\m}$ is not a
pole of $\zeta(t \la \m', \beta \ra )^{\gamma(\m')}$ \\
(since
the only possible pole of this function is $\frac{1}{\la \m', \beta \ra}\in \R$ and $t_\m \not \in \R$);
\item 
for all $\m' \in {\it  Ex }_-$, $t_{\m}$ is not a pole of 
$\zeta(t \la \m', \beta \ra )^{\gamma(\m')}$. \\
( if this were false, then   $\rho : = t_\m \la \m', \beta \ra$ would be a zero of $\zeta (s)$ satisfying:
$$ |\rho| = |t_\m| \cdot \la \m', \beta \ra 
=\frac{|\rho_0| \cdot \la \m', \beta \ra }{\la \m, \beta \ra}
\leq \frac{|\rho_0| \cdot |\m'| \cdot  (\sup_i\, \beta_i) }{|\m| \cdot 
  (\inf_i\, \beta_i)}\leq \frac{ a \cdot B}{2 \cdot |\m|}\leq \frac{ a }{2},$$
which is impossible);
\ee
By (\ref{gdeltam}) and the fact that $|\m| \le N_{r(\m)},$  we conclude that for each  
$\m \in {\it  Ex }_+ $ satisfying $|\m| \geq
B$, $t_\m$ is a zero of 
$Z_\beta (t).$ Since  $t_\m \rightarrow 0$ when $|\m| \rightarrow
+\infty,$ it follows that $\ds \{t_\m\}_{\{|\m| \geq B\}}$ contains  a
sequence of zeroes
of $Z_\beta (t)$ with accumulation point in 
$U(\delta_1)\,$ {\it if } $Z(h;\s)$ could be meromorphically extended to $V(h; \delta_0).$  As in the first subcase, this is
not possible.

\vskip .2 in

{\bf {{Case 2: ${\it  Ex }$ is finite}}} \par

\vskip .1 in

Set $\ds G(\s):=\bigg(\prod_{m\in {\it  Ex }} 
\zeta(\la \m, \s \ra )^{-\gamma(\m)}\bigg) Z(h;\s)$. We will prove that $G (\s) \equiv 1$ . \\

By choosing $r$ sufficiently large in the equation (\ref{gdeltaa}), we deduce that:
\be
\item $G(\s)$ is an Euler
product of the form $G(\s)= \prod_{p} \left(\sum_{\alpha \in \N_0^n} 
\frac{m_\alpha }{p^{\la \alpha ,\s \ra}}\right),$ where 
    $m_{\bold 0} = 1,$  and there exist $C,D>0$ such that \ 
   $m_\alpha  \leq C (1+ |\alpha|^D)$ \ for all $\alpha.$ 
\item $G(\s)$ converges absolutely in $V(h;0)=\cup_{r} V(h;\frac{1}{r})$.
\ee
Suppose that $G(\s) \not \equiv 1$. Then there exists 
$\alpha \neq \bold 0$ such that $m_\alpha \neq 0$. 
Now fix $\beta =(\beta_1,\dots,\beta_n) \in (0, \infty)^n$ as in Case 1. 
It follows that the   Euler product  
$$t\mapsto R_\beta (t):=G(t \beta)= 
\prod_{p} \bigg(\sum_{\alpha \in \N_0^n}\frac{m_\alpha }{p^{t\la
      \alpha ,\beta \ra}}\bigg)$$ converges absolutely in the halfplane
$\{ t\in \C : \Re(t)>0\}$.\par 
Set ${\cal  S}:=\{\alpha \in \N_0^n : m_\alpha \neq 0\}.$  
Since $\la \alpha , \beta \ra \rightarrow +\infty$ as 
$|\alpha|\rightarrow +\infty$, it is  clear that there exists  
$\nu \neq \bold 0 \in {\cal S}$ such that 
$\la \nu , \beta \ra 
=\inf_{\alpha \neq \bold 0 \in \cal S } \ \la \alpha , \beta \ra > 0$. We
fix this $\nu$ in the sequel.\\

\vskip .1 in

Let $N=\left[\frac{8 \la \nu , \beta \ra}{\inf_i \beta_i}\right]+ |\nu| + 1 \in \N$. Then we have for 
$\Re(t)> \frac{1}{2\la \nu ,\beta \ra}$ and uniformly in $p$:
\begin{eqnarray*}
\sum_{|\alpha| \geq N+1} \left|\frac{m_\alpha }{p^{t\la \alpha ,
      \beta\ra}}
\right|
&\ll  & 
\sum_{|\alpha| \geq N+1} \frac{|\alpha |^D}
{p^{\Re(t) \cdot |\alpha| \cdot (\inf_i\, \beta_i)}}\\
&\ll  & 
\sum_{|\alpha| \geq N+1} \frac{|\alpha |^D}
{p^{\Re(t) \cdot \frac{|\alpha|}{2} \cdot (\inf_i \beta_i)}} \cdot 
\frac{1}
{p^{\Re(t) \cdot \frac{N+1}{2} \cdot (\inf_i \beta_i)}}\\
&\ll  &
\frac{1}{p^{\Re(t) \cdot \frac{N+1}{2}\cdot (\inf_i \beta_i)}}
\sum_{|\alpha| \geq N+1} \frac{|\alpha |^D}
{2^{ |\alpha| \inf_i \beta_i / 4  \la \nu ,\beta \ra }}\\
&\ll  &
\frac{1}{p^{\Re(t) \cdot \frac{N+1}{2} \cdot (\inf_i \beta_i)}} \ll  \frac{1}{p^2}.
\end{eqnarray*}

From this we deduce that 
$$R_\beta (t)=G(t \beta)= 
\prod_{p} \bigg(\sum_{\alpha \in \N_0^n\atop |\alpha| \leq N}
\frac{m_\alpha }{p^{t\la \alpha ,\beta \ra}}+V_N(p;t)\bigg),$$
where $t\mapsto V_N (p;t)$ is a holomorphic function that satisfies the bound  
$V_N(p;t) \ll_N p^{-2}$ uniformly in $p$ and all  
$t \in \C$ such that  
$\Re(t) >\frac{1}{2\la \nu ,\beta\ra}\,.$ Since this Euler product   
converges absolutely for $t=\frac{1}{\la \nu, \beta \ra}>0,$ it follows that 
$$\prod_p \bigg(1 + \sum_{0 < |\alpha| \le N} \frac{m_\alpha}{p^{t\la \alpha ,\beta \ra}}
\bigg)$$
 also converges absolutely  for $t=\frac{1}{\la \nu, \beta \ra}.$ However, since $|\nu| \le N$ it follows that 
$\sum_{p} \frac{ m_\nu  }{p^{t\la \nu ,\beta \ra}}\bigg|_{t=1/\la \nu, \beta\ra}$ must also converge, which
is {\it not } possible. Thus, we conclude that $G(\s)\equiv 1$.\\ 

As a result, we must have the following equation  for all $\s \in V(h;A)$: 
\begin{eqnarray*}
Z(h;\s)&=&\prod_{m\in {\it  Ex }}\zeta(\la \m, \s \ra )^{\gamma(\m)}
=\prod_{m\in {\it  Ex }}
\prod_{p}\left(1-p^{-\la \m, \s \ra}\right)^{-\gamma(\m)}\\
&=& \prod_{p}\prod_{m\in {\it  Ex }}\left(1-p^{-\la \m, \s
  \ra}\right)^{-\gamma(\m)}
= \prod_{p} h^*(p^{-s_1},\dots,p^{-s_n}),
\end{eqnarray*}
where $h^*(X)=h^*(X_1,\dots,X_n)=
\prod_{m\in {\it  Ex }}\left(1-X^\m \right)^{-\gamma(\m)}=
\prod_{m\in {\it  Ex }}\left(1-X_1^{m_1}\dots X_n^{m_n} \right)^{-\gamma(\m)}$.
Since the Euler product factorization is unique,  we conclude that 
$h(X)=h^*(X),$ which completes the proof of Theorem \ref{ester+}.
\hfill $\square$ 

When $h  = 1 + \sum_{\alpha \neq \bold 0} a_\alpha X_1^{\alpha_1}\dots X_n^{\alpha_n}\in
\Z[X_1,\dots,X_n],$  we also have the analog of  Corollary \ref{zehf}, whose notation is used below.
\begin{Corollaire}\label{ester+f}
The Euler product  $Z(h;\s)$ can be continued from $V^\#(h; 1)$ to   $\C^n$ as a meromorphic function 
if and only if $h$ is `cyclotomic'. In all other cases $V^\#(h;0)$ is a natural boundary  (that is, $Z(h;\s)$ can not be
continued meromorphically to $V^\#(h;\delta)$ for any $\delta < 0$).
\end{Corollaire}
{\bf Proof:} The hard part is to prove the necessity, that is, $h$ must be cyclotomic if a  
meromorphic extension to
$\C^n$ exists. As with Theorem 2,   we will show this even if there exists an extension into $V^\# (h; \delta)$ for some
$\delta < 0.$ By a permutation of coordinates, one can suppose that:\\ 
$\big\{k\in \{1,\dots, n\} : \exists a \in \N 
{\mbox { s.t. }} 
a \e_k \in S^*(h)  \big\} = \{1,\dots,r\}.$\\
If the set is empty, then $r = 0.$ 

Assuming the set is nonempty, define $c_1,\dots,c_r \in \N$  by setting $c_k = \inf \{ c >0 : c \e_k \in S^*(h)\},$   for each
$k=1,\dots,r. $ It is clear that    
  $c_k \e_k \in S^*(h)$ for each $1 \le k \le r.$
If $r=n$, then  $c_k \e_k \in
Ext\left(S^*(h)\right)$ for all $k.$ Setting, for any $\delta \in \R,$ \  $\delta' = \frac{\delta}{\max_k c_k}$ and $\delta''
=
\frac{\delta}{\min_k c_k},$ this implies   that  
$V(h;\delta'') \subset V^\#(h;\delta) \subset V(h;\delta')$ if $\delta \ge 0,$ while $V(h;\delta') \subset V^\#(h;\delta)
\subset V(h;\delta'')$ if $\delta < 0.$  The assertion in Corollary \ref{ester+f} therefore follows immediately from the proof
of  Theorem \ref{ester+}.\par
Let us then suppose that $r<n$. 
We set
\begin{eqnarray*}
h^*(X_1,\dots,X_n)&:=& h(X_1,\dots,X_n) \prod_{k=r+1}^n
(1-X_k)\\
&=& 
\bigg(1+\sum_{\alpha \in S^*(h)}
a_\alpha X^{\alpha}\bigg) 
\bigg(\sum_{\eps \in \{0,1\}^{n-r}}
(-1)^{|\eps|} \prod_{k=r+1}^n X_k^{\eps_k}\bigg) \qquad (\eps = (\eps_{r+1},\dots, \eps_n))\\
&=&
1+\sum_{\alpha \in S^*(h)} a_\alpha  X^\alpha -\sum_{k=r+1}^n X_k \\
& & +\sum_{\eps \in \{0,1\}^{n-r}\atop |\eps |\geq 2}
(-1)^{|\eps|} \prod_{k=r+1}^n X_k^{\eps_k} +
\sum_{\alpha \in S^*(h)} \sum_{\eps \in \{0,1\}^{n-r}\atop |\eps |\geq 1}
(-1)^{|\eps|} a_\alpha X^\alpha \prod_{k=r+1}^n X_k^{\eps_k}.
\end{eqnarray*}
For each $ k  \ge r+1,$ set $c_k=1$.\\
It is then clear that    $c_k \e_k \in S^*(h^*)$ for all $k = 1,\dots,
n.$ Moreover, it follows immediately that 
$\sigma_k> \frac{1}{c_k}$
for each
$k \ge 1$ implies:
\begin{equation}\label{zhzh*}
 Z(h;\s) \prod_{k=r+1}^n {\zeta(s_k)}^{-1}= Z(h^*;\s).
\end{equation}

Suppose that there exists $\delta_0 < 0$ such that $ \s \mapsto
Z(h;\s)$ can be
meromorphically continued to  $V^\#(h;\delta_0).$
We set 
$\ds \delta_1=\frac{\delta_0}{2} \bigg(\sup_{\alpha \in S^*(h)} 
\bigg(\sum_{k=1}^n \frac{\alpha_k}{c_k}\bigg)\bigg)^{-1}<0$.
It is easy to check (exercise left to reader) that $V^\#(h^*;\delta_1) \subset
V^\#(h;\delta_0)$. This together with the relation (\ref{zhzh*}) 
then implies that $\s\mapsto Z(h^*;\s)$  can be
meromorphically continued to $V^\#(h^*;\delta_1)$.
Since  there
exists, for each $k,$ an integer $c_k \ge 1$ such that 
$c_k \e_k \in S^*(h^*),$  the proof  in the case   $r=n$ applies, from which it follows that   $h^*$ is a cyclotomic
polynomial.  The definition of $h^*$ then implies that the polynomial
 $h$ is also cyclotomic.  This completes the proof of Corollary \ref{ester+f}. $\square$\\

{\bf Remark:} \ Thus, for $h$ as above and  not cyclotomic,  the position of $\partial V(h;0)$ for a polynomial
can  differ rather significantly from that for an analytic function. Indeed, for the latter, $\partial V(h;0)$ is
always a union of coordinate hyperplanes, whereas for the former, $\partial V(h;0)$ need {\it not} be a subset of
$\partial [0,
\infty)^n.$ The situation is much less clear when 
$h = h (X_1,\dots, X_n, X_{n+1})$ and
$Z(\s) =
\prod_p h (p^{-s_1},\dots, p^{-s_n}, p)$ is defined as in Theorem 1 (see \S 3.2).  

 \section{An application  in diophantine geometry}

In    problems with a multiplicative structure, one often 
wants to estimate a counting function 
associated to a  multiplicative function. Since our interest is multivariate in nature, 
  a  multiplicative function for us refers to   any function
$F:\N^n \longrightarrow \C$  such that if  $(m_1,\dots,m_n),$ $(m_1',\dots,m_n')\in \N^{n}$ satisfy 
$gcd\left( {\it lcm} (m_1,\dots,m_n), {\it lcm} (m_1',\dots,m_n')\right)$ $=1,$ then  we have:
$$F(m_1 m_1',\dots,m_n m_n')=F(m_1,\dots,m_n)F(m_1',\dots,m_n').$$
To extract information about the averages of   $F(m_1,\dots, m_n)$ when the vectors $(m_1,\dots, m_n)$ are confined to 
a family of increasing   sets, 
it is often useful to study the analytic properties of an associated multivariate  Dirichlet series whose coefficients
are the values of the multiplicative function. Typically,  multiplicativity  implies that  such a Dirichlet series will be 
expressible as an absolutely convergent  Euler product in some domain. When the values of $F$ are integral, one might expect 
that the results of \S 1 could then be applied. 

We show that this is indeed possible   with an   example   from toric  geometry.  
In this case,    the multiplicative function that appears quite naturally  is not only  integral valued but also 
satisfies a  special invariant property. \S 2.1 shows how Theorem 1 or Corollary 1 can be used to deduce pertinent properties
of   Dirichlet series  whose coefficients are determined by an ``invariant  multiplicative function". \S 2.2 defines the
invariant multiplicative function that can be associated to any toric variety, and discusses some necessary ingredients
for the
main result  of this section that is the subject of  \S 2.3. In \S 2.4, we then apply this result to give the explicit 
asymptotic for a general problem in multiplicative number theory that is equivalent to a counting problem in toric
geometry. 
 
\subsection{Properties of a   Dirichlet series with multiplicative coefficients}

\begin{Definition} \label {def1} An integral multiplicative function $F:\N^n \to \Z$ is said to be invariant if 
 the function $\nug = (\nu_1,\dots,\nu_n) \in \N_0^n \to F(p^{\nu_1},\dots, p^{\nu_n})$ is independent of $p.$
In this event,   the function $f(\nug): = F(p^{\nu_1},\dots, p^{\nu_n})$ is well defined, and  called the index function of
$F.$

The support of   $f$ is the set \ $S(f):= \{\nug \in \N_0^n \mid f(\nug) \neq 0 \}.$ The index of $S^*(f)$  (see Notations
part (8)) is denoted 
$\iota (f):=\iota (S^*(f)).$  For each $\alpha \in R(S^*(f)),$   we define:
\be 
\item $K_+(f;\alpha):=\{\nug \in K(S^*(f);\alpha) : \ 
f(\nug)>0 \}$;  
\item $K_-(f;\alpha):=\{\nug \in K(S^*(f);\alpha) : \ f(\nug) < 0 \}$;
\item $J(\alpha)=\{\e_j :\  \alpha_j =0\}$.
\ee 
\end{Definition}

We will also {\it assume}  that the index function $f$ satisfies this   property:
\begin{equation}\label{9}
{\mbox {\it if   $S^*(f)$ is  not finite, then  $S^*(f)\cap (\R \e_i) \neq \emptyset$ for each $i = 1,\dots, n.$}}
\end{equation}

By applying Theorem 1 we   prove the following. 
\begin{Theoreme}\label{application1}
Let $f$ be the index function of an invariant integral multiplicative function $F.$ In addition to (\ref{9}), we assume
that there exist  constants $C,D>0$  such that   $|f(\nug)|\leq C (1+| \nug |)^D.$  If there exists $\alpha \in
R(S^*(f))$ such that 
$K_-(f;\alpha)=\emptyset$, then there exist $\theta > 0$ and a polynomial 
 $Q \in \R[X]$  such that: 
 \begin{equation} \label{estimemultiple}
\sum_{(m_1,\dots,m_n)\in \N^{n} \atop \max_{i} m_i \leq t} 
F(m_1,\dots,m_n) = t^{\iota (f)} Q(\log t) + O(t^{\iota (f)-\theta})
\quad as ~~ t\rightarrow \infty,
\end{equation}
where the degree of $Q$ is at most 
$$q = \min_{\alpha \in R(S^*(f))} \bigg\{\#J(\alpha) + \sum_{\nug \in
K_+(f;\alpha)} f(\nug) 
 -rank \left\{\nug : \nug \in
  K_+(f;\alpha)\cup J(\alpha)\right\} \bigg\}.$$  
If, however, $K_-(f;\alpha) \neq \emptyset$ for each
$\alpha \in R(S^*(f)),$ then   (\ref{estimemultiple}) continues to hold   if we assume the Riemann Hypothesis.
\end{Theoreme}
{\bf Proof:} \ First we define the   zeta function: 
$$Z(F;\s):= \sum_{(m_1,\dots,m_n)\in \N^{n}}
\frac{F(m_1,\dots, m_n)}{m_1^{s_1}\dots m_n^{s_n}}.$$
The multiplicativity of $F$ and the bound on $f(\nug)$ imply that there exists a constant $D>0$
such that for all  $(m_1,\dots, m_n) \in \N^{n}$: 
$$F(m_1,\dots, m_n)\ll \prod_{p | m_1\dots m_n}
\left(1+v_p(m_1 \dots m_n)\right)^D \leq  {\tau (m_1\dots m_n)}^D,$$
where $\tau$ denotes the usual divisor function and $v_p$ denotes $p-$adic order. 
By using the standard bound $\tau (d)
\ll_\varepsilon d^\varepsilon\,, $  it follows
that    
$$F(m_1,\dots,m_n)\ll_\varepsilon {(m_1\dots m_n)}^\varepsilon\,.$$ 
Thus, $Z(F;\s)$ converges absolutely on  
$\Omega :=\{\s\in \C^n : \sigma_i >1 \ \ \forall i=1,..,n \}$.

Multiplicativity 
of $F$ then implies the following formula for all $\s \in \Omega:$
$$Z(F;\s)=\prod_{p}
\bigg(\sum_{\nug \in \N_0^n} \frac{F(p^{\nu_1},\dots,p^{\nu_n})}{p^{\la
        \nug ,\s \ra}}\bigg)=\prod_{p}
\bigg(1+\sum_{\nug \in \N_0^n\setminus \{0\}} \frac{f(\nug)}{p^{\la
        \nug ,\s \ra}}\bigg).$$
Defining  $h_f(X_1,\dots,X_n) = \sum_{\nug \in \N_0^n} 
f(\nug) X_1^{\nu_1}\dots X_n^{\nu_n},$ it is then clear that ${\cal E} (h_f) = {\cal E}(S^*(f)),$ and $Z(F;\s) = Z(h_f;\s)$
for $\s \in \Omega.$  Note that $h_f = h_{f,0}$ in the notation prior to Theorem 1. It  will also be convenient to abuse
notation by writing below  $V(f;\delta)$ for the sets denoted $V(h;\delta),$  with $h = h_f$, in \S 1.1.   

We can then apply Theorem 1 (and, in addition, the Remark that immediately follows its proof) if $S^*(f)$ is infinite,  or
Corollary 1 if
$S^*(f)$ is finite, to conclude the following.
\be
\item $\s\longrightarrow Z(F;\s)$ converges absolutely in the
  domain $V(f;1)$ and admits a meromorphic continuation to the set
  $V(f;0)$;
\item there exists  $\delta \in (0,1)$ and a  holomorphic bounded
  function  $G_\delta$ in $V(f;\delta)$  such that  
$$Z(F;\s)=  \prod_{\nug \in {\widetilde { S^*(f) }}} 
\zeta(\la\nug,\s\ra)^{f(\nug)}\, \cdot \, G_\delta(\s) \quad \forall \s\in V(f;1).$$
\ee

Let   $\alpha \in R(S^*(f))$ be given,  
and assume $\s\in \C^n$ satisfies $\sigma_i  > \alpha_i$ for each $ i=1,\dots,n.$ Then $Z(F;\s)$ converges
absolutely since the inequality 
$ \la \sigma ,\nu \ra > \la \alpha ,\nu \ra  \geq 1$ for all $\nu \in {\widetilde {S^*(f)}}$ follows from the 
assumption on $\s.$  

 The  argument to follow now has two parts. The first gives an explicit expression for the meromorphic continuation of 
$Z(F;\s)$  in a neighborhood of $\alpha.$ The second gives an 
explicit description of a  divisor containing $\alpha$ that could contain   components of the polar divisor of the
meromorphic continuation.

For each $\alpha \in R(S^*(f)),$ we first define:
$${\cal H}_\alpha(\s):=Z(F;\alpha +\s)~\prod_{\nug \in K_+(f;\alpha)} \la\nug,\s\ra^{f(\nug)}. $$ 
Then for all $\s \in V(f; 0)$
\begin{eqnarray*}
{\cal H}_\alpha(\s)&=&  \prod_{\nug \in K_+(f;\alpha)} 
\left(\la\nug,\s\ra \zeta(1+\la\nug,\s\ra)\right)^{f(\nug)}\, \cdot \,
 \prod_{\nug \in K_-(f;\alpha)} 
 \zeta(1+\la\nug,\s\ra)^{f(\nug)}\\
&& \times 
  \prod_{\nug \in {\widetilde {S^*(f)}} \setminus K(S^*(f);\alpha)} 
\zeta(\la\nug,\alpha \ra +\la\nug,\s\ra)^{f(\nug)}\, \cdot \,
G_\delta(\alpha +\s).
\end{eqnarray*}
 Using  classical properties of the Riemann zeta function and
assuming the Riemann hypothesis in the case $K_-(f;\alpha)\neq \emptyset$, 
it is easy to
see that there exist $\delta_1,\delta_2 >0$ such that:
\be
\item $\s\mapsto {\cal H}_\alpha(\s)$ is holomorphic in the set 
$$V(f;-\delta_1):=
\{\s\in
\C^n :   \la \sigma ,\nug\ra > -\delta_1 \ \ \forall \nug \in {\widetilde {S^*(f)}}\}.$$
In particular,  $G_\delta(\alpha + \s)$ is holomorphic in some $V(f,-\delta_1)$ if $S^*(f)$ is infinite since the
hypothesis (\ref{9}) implies each $\alpha_i  > 0.$ If $S^*(f)$ is finite, then this property holds, even if some $\alpha_i =
0,$ by the proof of Corollary 1. 
\item \ \ $ {\cal H}_\alpha(\s)\ll_\epsilon \prod_{\nug \in K^+(f;\alpha)} 
\left(1+| \la \nug, \tau \ra|\right)^{f(\nug)(1-\delta_2\min\{0, \la\nug,\sigma \ra)\})}
~\bigg(1+(\sum_{i=1}^n |\tau_i|)^\epsilon\bigg),$ \newline
where the implied constant is independent of  $\s \in V(f;-\delta_1)$. 
\ee
This gives an explicit expression for the meromorphic continuation of $Z(F;\s)$ in a neighborhood of $\alpha$ as
\begin{equation} \label{merocont}
Z(\alpha + \s) = \frac{{\cal H}_\alpha (\s)}{\prod_{\nug \in K_+(f;\alpha)} \la\nug,\s\ra^{f(\nug)}}\,.\end{equation}
The second part of the argument now follows easily. This equation also shows that the divisor  
 $${\cal D_\alpha}: = \sum_{\nug
\in K^+(f;\alpha)} f(\nug) \cdot \bigg\{ \alpha + \{\la \nug, \s \ra = 0\} \bigg\}$$   
could contain components of the polar divisor of the quotient. 

The growth estimate  in part 2   says that the quotient in (\ref{merocont}) 
grows  at a
polynomial rate in $\tau$ when $\sigma$ is confined to any bounded neighborhood of $\bold 0$ that lies inside 
$V(f;-\delta_1)_{\R}\,,$ {\it and} \ $\s$ remains a positive distance away ${\cal D_\alpha}.$  

{\bf Remark:}\ It is also important to observe that the preceding argument can be easily extended to any  point in
${\widetilde {S^*(f)^o}}.$ Since this set is   convex, its boundary can be thought of as the set of vectors $\xi$ such
that
$\la \xi, \nu \ra = 1$ for some $\nu \in S^*(f).$ The set $K(S^*(f), \xi)$ is then seen to equal the support plane
to $\partial {\widetilde {S^*(f)^o}}$ in the direction of $\xi.$ 
One can therefore think of this boundary as a first approximation to 
the Newton polyhedron of the polar
divisor of $Z(F;\s)$ in the sense of  (\cite{lichtinduke1},7.1).  Its support plane
$\ell$ in the particular direction of the diagonal   $(1,\dots, 1)$ intersects the boundary exactly in  the set
$R(S^*(f)).$\hfill $\square$

By an iteration of Perron's lemma in $\C^n,$ it follows that for any $t \notin \N$ and $c \gg 1,$  
\begin{eqnarray}\label{perron}
\frac 1{(2 \pi i)^n}
\int_{\{\sigma = (c,\dots, c)\}} Z(F;\s)  t^{s_1+\cdots +s_n} \  \frac{d s_1 \cdots d s_n} {s_1 \cdots s_n} &=&
\sum_{(m_1,\dots,m_n)\in
\N^{n}
\atop
  m_i \leq t \ \forall i} F(m_1,\dots, m_n) \\
&=& \sum_{(m_1,\dots,m_n)\in \N^{n} \atop \max_{i} m_i \leq t} F(m_1,\dots, m_n)
 \nonumber\,.\end{eqnarray} 
(If $t \in \N,$ then one needs to multiply $F(m_1,\dots, m_n)$ by $1/2$ if $F(m_1,\dots,m_n) = t.$)  
Applying the preceding Remark, we can then use the method described in [ibid., section 7] (also see 
\cite{lichtinduke2}, Appendix) to deduce the asymptotic behavior in
$t$ of the average on the right side of (\ref{perron}). We do so by replacing the Newton polyhedron of the polar
divisor of $Z(F;\s)/s_1\cdots s_n$ by $\partial \tilde X^o$ where $X = S^*(f) \cup \{\bold e_j\}_1^n.$ The dominant term,
{\it if it is nonzero} (which need not occur, {\it even if} this set equals the Newton polyhedron of the
polar divisor!), will be given by
$t^{|\xi|} Q(\log t),$ 
where $Q$ is a polynomial of degree at most the integer $q$ defined in the statement of the Theorem, and $\xi$ is any
vertex of $\partial \tilde X^o$ that also belongs to the    support plane
$\ell$ defined above. {\it By  definition}, it   follows that $|\xi| = \iota (f).$ This completes the proof of the
Theorem. 

{\bf Remark:} For the reader
more comfortable with purely 
$1-$ variable methods, it is worthwhile to indicate that once one knows (\ref{merocont}),  
\cite{bretechecompo} has derived the same asymptotic in (\ref{estimemultiple}) by using  a  procedure  of  standard Cauchy
residue techniques suitably iterated. An advantage of  that method is that it also gives a condition that is easy to
check, and sufficient to show that the polynomial $Q (\log t)$   in (\ref{estimemultiple}) is {\it nonzero} (see \S 2.3).
On the other hand, the method  of
\cite{lichtinduke1} can  be used to derive the expected dominant asymptotic in $n$ independent parameters
$t_1,\dots, t_n$ in place of a single parameter $t.$ We will not, however, develop this point here. 

\subsection{An invariant multiplicative function associated to a projective toric variety}

In this and the next subsection,   $\bold A$ denotes a  $d\times n$   matrix  
with entries  in $\Z,$   whose rows   
$\a_j=(a_{j,1},\dots,a_{j,n})$ each satisfy the property that 
$\sum_{i} a_{j,i} = 0.$  One can then define the following objects:  
\begin{eqnarray*}
&{}& {\mbox {rational points of a projective toric variety}} \ X(\bold A)\\
&{}& := \{(x_1:\dots:x_n) \in \P^{n-1}(\Q) :   
\prod_{\{i: a_{j,i} \ge 0\}} \, x_i^{a_{j,i}} = \prod_{\{i: a_{j,i} < 0\}} \, x_i^{-a_{j,i}}\ \ \forall j\};\\
&{}& {\mbox {an  open subset}} \ U(\bold A)\\
&{}& :=\{(x_1:\dots:x_n) \in X(\bold A) : x_1\dots x_n \neq 0 \};\\
&{}&\\
&{}&{\mbox { a subset of $\ker \bold A,$}}\quad
 T(\bold A) \\
&{}& := \{\bold \nu \in \N_0^n: \bold A(\bold \nu) = \bold 0 \ \ {\mbox { and }}\ \  \prod_i \nu_i =
0\}.\end{eqnarray*}
Following the  idea in \cite{salberger},  we define a function $F_{\bold A}:\N^n \to \Z$  by setting:
\be
\item
 $F_{\bold A} (m_1,\dots,m_n)=1$ \ if \ $\gcd(m_1,..,m_n)=1$ \  and \   
$\prod_i m_i^{a_{j,i}} = 1$  \ $\forall j\le d$,
\item $F_{\bold A}(m_1,\dots,m_n)=0$ \ if not.
\ee
It is clear that $F_{\bold A}$ is multiplicative, $F_{\bold A}(m_1,\dots, m_n) = 1$ implies $(m_1:\dots:m_n) \in
U(\bold A),$ and that for each $p$ and all $\bold \nu \in \N_0^n$,
$$F_{\bold A} (p^{\nu_1},\dots,p^{\nu_n})= 1 \ {\mbox { iff }} \ \bold \nu \in T(\bold A).$$
Thus, $F_{\bold A}$ is   invariant, and its index function is the characteristic function
of
$T(\bold A).$

By the multiplicativity of $F_{\bold A},$ we obtain, exactly as in \S 2.1, that for all  $\s \in \Omega,$  
$$Z(F_{\bold A};\s): =\sum_{(m_1,\dots,m_n)\in
  \N^n} \frac{F_{\bold A}(m_1,\dots,m_n)}{m_1^{s_1}\dots m_n^{s_n}}= 
\prod_{p} \bigg(\sum_{\bold \nu \in T(\bold A)} 
\frac{1}{p^{\la \bold \nu, \s \ra}}\bigg).$$
{\bf Note:} The relation between $Z(F_{\bold A};\s)$ and a ``generalized" height zeta function on $U(\bold A)$ is explained in
\S 2.3 (see (\ref{zetarel})).\hfill $\square$

Set
\begin{equation}\label{10}
 h_{\bold A} (X): = \sum_{\bold \nu \in T(\bold A)} X^{\bold \nu}   
\end{equation}
to denote the function whose coefficients are determined by the index function of $F_{\bold A}.$ The only  thing that 
we know of for sure about $h_{\bold A}$  is that it is  analytic  on the unit polydisc
$P(1)$ in $\C^n.$ This, however, does not even allow us to apply Theorem 3.  It will therefore be necessary to
understand this function   much more precisely.  
 
The crucial property   is the following.
\begin{Definition}
An analytic function $h$   on   $P(1)$ is unitary if   
there exist a finite set  $K \subset \N_0^n \setminus\{0\},$   positive integers $\{c(\bold \nu)\}_{\bold \nu\in K},$ 
 and a polynomial $W \in \Z[X_1,\dots,X_n]$, 
such that for all $X\in P(1)$:
$$h(X) = \bigg(\prod_{\bold \nu \in K}
  (1-X^{\bold \nu})^{-c(\bold \nu)}\bigg) W(X).$$
The data 
$\left(K;\la c(\bold \nu)\ra_{\bold \nu\in K};W\right)$ determines a presentation of $h$  when
$1-X^{\bold \nu}$ does not  divide  $W(X)$   for each $\bold \nu \in K.$ 
\end{Definition} 

The particular result we will need in \S 2.3 will then be as follows.
\begin{Lemme}\label{lemdefbis} The function $h_{\bold A},$ defined in (\ref{10}), is unitary.
\end{Lemme}
 
Lemma {\ref{lemdefbis}} is a simple consequence of a   more general result which analyzes the behavior of an analytic 
function, all of whose monomial exponents   belong  to an  affine plane 
$$T(\bold A, \bold b): = \{\bold \nu \in \N_0^n: \ \bold A (\bold \nu) = \bold b \}.$$
\begin{Lemme}\label{lemdef}
For any integral \ $d\times n$ matrix $\bold A$ (the rows of which  need not  sum to $0$!), and any $\b \in
\Z^d,$ the function 
$$h_{\bold A;\b}(X):=\sum_{\bold \nu \in T(\bold A;\b)} X^{\bold \nu}$$
 is unitary.
\end{Lemme}

{\bf Proof that Lemma \ref{lemdef} implies Lemma \ref{lemdefbis}:}\par
For all $X=(X_1,\dots,X_n) \in P(1)$ we have:
\begin{eqnarray*}
h_{\bold A}(X)
&=&   \sum_{\bold \nu \in T(\bold A; \bold 0)\atop \nu_1\dots\nu_n =0} 
X^{\bold \nu} =
 \sum_{\bold \nu \in T(\bold A;\bold 0)} 
X^{\bold \nu} - \sum_{\bold \nu \in T(\bold A;\bold 0)\atop \nu_1\geq
  1,\dots,\nu_n \geq 1} 
X^{\bold \nu}\\
&=& (1-X_1\dots X_n) h_{\bold A, \bold 0}(X).
\end{eqnarray*}
Since Lemma \ref{lemdef} says that $h_{\bold A, \bold 0}$ is unitary
it follows   that $h_{\bold A}$ is also unitary. \hfill $\square$

{\bf Proof of Lemma \ref{lemdef}:}\par
We shall prove the lemma by induction on $n$.\par
For $n=1$ the result is trivially true.\par 
Let $n\geq 2.$ The induction hypothesis allows us to assume   that for any $m < n,$ any  $d \times m$ integral matrix $\bold
A',$   and any $\bold b' \in \Z^d,$ we have that $h_{\bold A', \bold b'}(X_1,\dots, X_m)$ is unitary.

Now, let $\bold A$ be a $d \times n$ integral matrix,  and  $\b=(b_1,\dots,b_d)\in \Z^d.$ It   suffices to assume
that  $T(\bold A, \bold b) \neq \emptyset$ since the proof of Lemma 2 is trivial when $T(\bold A, \bold b) = \emptyset.$
  
It will be convenient to distinguish two cases:

{\bf Case 1 :  \ $\{\bold 0\} \subsetneq  T(\bold A;\bold 0).$} 

We choose and fix   $\bold \alpha \neq \bold 0 \in T(\bold A;\bold 0)$ in the following. For any $I\subset \{1,\dots,n\}$,
we define 
$$L(I,\bold \alpha):=\{\bold \nu \in T(\bold A, \bold b):  \nu_i \geq
  \alpha_i \ {\mbox {\it iff }}\ i \in I\},$$
and 
\begin{equation}\label{fix}
h_{\bold A, \bold b}(I;\bold \alpha;X):= \sum_{\bold \nu \in L(I;\bold \alpha)} X^{\bold \nu}.
\end{equation} 
If $L(I,\bold \alpha) = \emptyset,$ the value is defined to be $0.$  A straightforward  calculation then shows:
\begin{equation}\label{decalage}
 (1-X^{\bold \alpha}) h_{\bold A, \bold b}(X)= \sum_{I\subset \{1,\dots,n\}
 \atop I \neq \{1,\dots,n\}} h_{\bold A,\bold b}(I;\bold \alpha;X) \quad \forall  X \in P(1). 
\end{equation}
So, we need to show that each $h_{\bold A, \bold b}(I;\bold \alpha;X)$  is unitary.
By permuting coordinates, it suffices to prove this for  any  $I_q:=\{ 1,2,\dots,q\}$ with $q\leq n-1$. 

To express the necessary equation in a concise manner, we first introduce the following notations:  
\be
\item $X=(Y,Z)$ with
$Y=(X_1,\dots,X_q)$ and  $Z=(X_{q+1},\dots,X_n);$ 
\item   $\x' = (x_1,\dots,x_q)$ and  $\x'' = (\x_{q+1},\dots, \x_n),$ for any $n-$vector $\x,$ and $\bold A'$ is the
$d \times q$ matrix with rows
$\bold a_j' = (a_{j,1},\dots, a_{j,q})$ for each $j \le d;$ 
\item ${\cal D} (= {\cal D(\bold \alpha)}) := \left\{\bold \nu''=(\nu_{q+1},\dots,\nu_n) 
\in  \prod_{i=q+1}^n \{0, 1, 2,\dots, \alpha_i -1\} \right\}$;
\item $\forall \bold \nu'' \in {\cal D}$,  
$$\l(\bold \nu''):=\bigg( b_1-\la  \a_1'', \bold \nu''\ra - \la \a_1',\bold \alpha' \ra,\dots, b_d-\la  \a_d'', \bold \nu''\ra
-
\la \a_d',\bold \alpha' \ra \bigg).$$  
\ee

We then observe  that  for all $X=(Y,Z) \in P(1)$,
\begin{eqnarray*}
h_{\bold A, \bold b}(I_q;\bold \alpha;X)
&=& \sum_{\nu \in T(\bold A, \bold b) \atop 
\forall i \le q \ \nu_i \geq \alpha_i  \ {\it and} \ 
\forall i > q \ \nu_i < \alpha_i} X^\nu \\
&=& \sum_{\nu''=(\nu_{q+1},\dots,\nu_n) 
\in  {\cal D}} 
\sum_{\mu=(\mu_1,\dots,\mu_q) \in \N_0^q \atop 
(\alpha'+\mu, \nu'') \in T(\bold A, \bold b)} Y^{\alpha'+\mu}
Z^{\nu''}\\
&=& \sum_{\nu''=(\nu_{q+1},\dots,\nu_n) 
\in  {\cal D}} Y^{\alpha'} Z^{\nu''}
\sum_{\mu=(\mu_1,\dots,\mu_q) \in  
T(\bold A', \l(\bold \nu''))} Y^{\mu}\,.
\end{eqnarray*}

So the following equation is true:
\begin{equation}\label{11} h_{\bold A, \bold b}(I_q;\bold \alpha;X)=\sum_{\nu'' \in {\cal D}} Y^{\bold \alpha'} Z^{\bold
\nu''}  h_{\bold A',{\l}(\nu'')} (Y). \end{equation}
We conclude by induction. 

{\bf Case 2 :  $T(\bold A;\bold 0)=\{\bold 0\}.$}

Since $T(\bold A, \bold b) \neq \emptyset,$  there exists $\bold \gamma \in T (\bold A;\bold b)$. We begin by observing
that:\\
 $\bold \nu \in T (\bold A;\b)$ 
is equivalent to one of the two following conditions :
\be
\item $\bold \nu =\bold \gamma$ (i.e. $\nu_i \ge \gamma_i$ \   $\forall i$ implies $\bold \nu - \gamma \in T(\bold A;\bold 0) =
\{\bold 0\}$);
\item $\bold \nu \in T(\bold A;\b)$ and 
$\exists i\in \{1,\dots,n\}$ such that $\nu_i < \gamma_i\,.$ 
\ee
This observation implies that for all $X\in P(1)$:
$$h_{\bold A,\b}(X)= X^{\bold \gamma}+\sum_{I\subset \{1,\dots,n\}
 \atop I \neq \{1,\dots,n\}} h_{\bold A,\b}(I;\bold \gamma;X)$$
where each  $h_{\bold A,\b}(I;\bold \gamma;X)$ is defined as in (\ref{fix}), replacing $\alpha$ by $\gamma.$
We now conclude by induction as in Case 1. 
This completes the proof of Lemma \ref{lemdef}.\hfill
$\square$  

{\bf Remark:} The proof of Lemma \ref{lemdef} actually gives an explicit procedure to find a presentation of $h_{\bold A}.$
This is useful to find  the polyhedron of $h_{\bold A}$ in specific examples, as \S 2.4   shows.   

\subsection {Analytic properties of a generalized height
zeta function for a toric variety}

We first recall that to any ``height"
vector, that is, any 
$\bold \beta = (\beta_1,\dots, \beta_n) \in (0, 1)^n$ satisfying
$\sum_i
\beta_i = 1,$ one can define  a height function $H_\beta$ on   $U(\bold A)$  by choosing the unique  
representative
$(x_1:\dots:x_n)$ of  a   point $\x \in U(\bold A),$ satisfying the properties that each 
$x_i \in \Z$ and  $gcd (x_1,\dots, x_n) = 1,$  and then setting
$$H_{\bold \beta} (\x) := \prod_i |x_i|^{\beta_i}.$$
The height zeta function on $U(\bold A)$ is a function of the complex variable $s$ and defined as the series
$$Z_{\bold \beta}(s) = \sum_{\x \in U(\bold A)} H_{\bold \beta} (\x)^{-s}.$$
Rather than focus upon a single choice of $\bold \beta,$ it is reasonable to look for   a single zeta function that 
contains  the information encoded by all the   $Z_{\bold \beta}.$ The natural choice is to define the ``generalized" height
function on $U(\bold A)$ by setting $\s = (s_1,\dots, s_n) \in \C^n$ and defining
$$H (\x, \s): = \prod_i |x_i|^{-s_i},$$
where the same choice of representative   for a point $\x$ is used as above. The corresponding generalized height zeta
function is then the multivariate Dirichlet series:
$$Z_{U(\bold A)}(\s) = \sum_{\x \in U(\bold A)} H(\x,\s).$$ 
It is clear that $Z_{U(\bold A)}$ is absolutely convergent on the open set $\Omega$ (see \S 2.1), and that $Z_{U(\bold
A)}(\bold
\beta s) = Z_{\bold \beta}(s)$ if $\sigma \gg_{\bold \beta} 1.$ 

Moreover, defining the constant
$$C(\bold A) := \frac 12 \cdot \# \bigg\{ \epsilon \in \{\pm 1\}^n : \prod_{i=1}^n  \epsilon_i^{a_{j,i}}=1 \ {\mbox { for
all }} j = 1,\dots, d \bigg\},$$
and recalling the definition of $F_{\bold A}$ from \S 2.2, it is easy to check that  $\s \in \Omega$   implies
\begin{equation} \label{zetarel}
Z_{U(\bold A)} (\s) = C(\bold A) \cdot Z(F_{\bold A};\s).\end{equation} 
Thus, the analytic properties of $Z(F_{\bold A};\s)$ are equivalent to   those of $Z_{U(\bold A)}(\s),$ and by specializing
$\s \to \bold \beta \cdot s$ we can infer   properties of each  $Z_{\bold \beta}(s).$  
 
The essential first step needed to deduce the analytic properties of $Z(F_{\bold A};\s)$ is given by Lemma \ref{lemdefbis}
in
\S 2.2. This insures that there is a presentation of $h_{\bold A}(X)$ as a rational function:
\begin{equation}\label{hA} h_{\bold A}(X) = \prod_{\bold \nu \in K} (1 - X^{\bold \nu})^{-c(\bold \nu)} \,\cdot \, W
(X).\end{equation}

{\bf Note.} \ Although $K$ and $W$ certainly depend upon $\bold A,$ the notation will not indicate this for the sake of
simplicity. The reader should not find this confusing.\hfill $\square$

Since both $h_{\bold A}(X)$ and each $(1 - X^{\bold \nu})^{-c(\bold \nu)}$ equal $1$ when $X = \bold 0,$ it is clear that
$W$ is a polynomial with integer coefficients that satisfies $W (\bold 0) = 1.$ Thus, Corollaries 1, 2 apply to the
Euler product
$Z(W;\s) = \prod_p W (p^{-s_1},\dots, p^{-s_n}).$

 For every $\delta \in \R$, define  
$\ds V(\delta):= 
\{\s \in \C^n \,: \, \la \bold \nu, \sigma \ra  > \delta  \ \ \forall \bold \nu \in S^*(W) \cup K \}.$
It is then clear that \ $Z_{\bold U(\bold A)} (\s)$ converges absolutely
in $V(1)$ and satisfies :
\begin{equation}\label {zetafunc} Z_{\bold U(\bold A)} (\s) = C(\bold A) \cdot 
\bigg(\prod_{\m\in K} \zeta\left(\la
  \bold \nu,\s\ra\right)^{c(\bold \nu)}\bigg) \cdot Z(W;\s).\end{equation}
Outside $V(1),$ Corollaries 1, 2 (whose notations are
used below) can now be immediately applied to tell us the following.
\begin{Theoreme}\label{application2} \hfill \newline 
\medskip
1.\ \ $\s \mapsto Z_{\bold U(\bold A)} (\s)$ can be
meromorphically continued to  $V(W ;0);$  

2.\ \ $\s \mapsto Z_{\bold U(\bold A)} (\s)$  can be
meromorphically continued to $\C^n$ if and only if  
$W $ is cyclotomic;

3.\ \ if $W $ is not cyclotomic, then $\partial V(W ;0)$ is
the  natural boundary of meromorphic continuation;

4.\ \ for any  height vector $\bold \beta,$     the  height zeta function
$Z_{\bold
\beta}(s)$ is either  meromorphic on
$\C,$ or, if not,  can  be meromorphically continued (at least) into the halfplane $\{\sigma > \eta_{\bold \beta}\},$
where
$\eta_{\bold
\beta}$ is   the point of intersection of the line $\{\bold \beta \cdot \sigma\}$ with 
$\partial V(W, \bold 0)_{\R}.$  
\end{Theoreme}

Manin's conjecture, applied to  a smooth toric model $X'(\bold A)$ of $X(\bold A),$   asserted a very precise asymptotic 
for the average   of the ``anticanonical height" function on $U(\bold A),$ when viewed as a dense   torus  on
$X'(\bold A).$ This is the   function, in down to earth terms, equal to
$$N_\infty (U(\bold A), t): = C(\bold A) \cdot \sum_{(m_1,\dots,m_n)\in \N^{n} \atop \max  m_j \leq t}
F_{\bold A}(m_1,\dots,m_n).$$   
The dominant term, as $t \to \infty,$ was conjectured to be asymptotically equivalent  to $C t \log^b t,$ where $C > 0$
had a specific expression as a product of certain volumes,  and  $b$ is one less than the rank of the Picard group of
$X(\bold A).$  The original conjecture was   proved by Batyrev-Tschinkel \cite{batyrev}. Salberger (\cite{salberger}, 11.1)
then used the theory of universal torsors to prove the asymptotic with an error term $O_{\epsilon}\big(t \log^{b - \frac 12 +
\epsilon}\,(t)\big),$ assuming the anticanonical bundle was ample. A bit later,    de la Bret\`eche
\cite{bretechetorique}   used   Salberger's work   and his own Tauberian theorem in
\cite{bretechecompo} to prove the   asymptotic   with a strictly smaller order (in the exponent for $t$) error term. 

We are able to extend this analysis by giving explicit asymptotics for many different kinds of   averages of $F_{\bold
A}(m_1,\dots, m_n).$ Given a   vector $\gamma = (\gamma_1,\dots, \gamma_n) \in (0, \infty)^n$ set
$$N_{\gamma} (U(\bold A), t) = C(\bold A) \sum_{1 \le m_i \le t^{\gamma_i} \atop  i = 1,\dots,n} F_{\bold
A}(m_1,\dots, m_n).$$
This     counts the points in $U(\bold A)$ in a family of boxes whose lengths in   different coordinate 
directions grow at different rates, according to the values of the components of $\gamma.$

The discussion to follow will prove an explicit (and nonzero) asymptotic for $N_{\gamma}(U(\bold A), t),$
whenever $\gamma$ is a ``generic" vector (see Theorem \ref{coefficient}). Because we have emphasized constructions
associated to a Newton polyhedron, it is natural that we  should express the dominant term   in   terms of a
polyhedron that is intrinsic to the problem. For our purposes, this   equals  the boundary of the dual of $K \cup
S^*(W).$ It is important to emphasize here that this can be computed {\it without} recourse to constructing an explicit
desingularized model of
$X(\bold A).$ Sometimes, at least, there are    computational advantages to this, as   \S 2.4 shows.

There are two parts to finding the asymptotic  of $N_{\gamma}(U(\bold A),t).$ The first part (see Theorem
\ref{application2bis}) proves  a necessary sharpening of Theorem 3. This shows that  the boundary of the dual of $K \cup
S^*(W)$ is {\it the} Newton polyhedron of the polar divisor of
$Z(U(\bold A);
\s)$ (in the sense of (\cite{lichtinduke1}, 7.1)),
{\it not merely} a first approximation. To prove this fundamental property, we exploit the fact   that  there is 
additional information built into the right sides of (\ref {hA}), (\ref{zetafunc}) than is available in general. The second
part (see Theorem
\ref{coefficient}) shows that the expected dominant term in the asymptotic is genuinely nonzero, and characterizes, as
well, the degree of the polynomial
$Q$ in polyhedral terms. For this, we use   some ideas from [ibid.], and a crucial nonvanishing property (see
(\ref{nonzero})) that is key to the proof of Theorem \ref{application2bis}. This then allows us to apply the Tauberian
theorem of de la Bret\`eche [op cit.]. 

{\bf Remark:} The reader should note that our results give considerably more information about the polar divisor of
$Z(U(\bold A);\s)$ than has been established in the preceding work cited above. In
particular, the earlier proofs of the asymptotic  have   all been based upon the ability to prove that 
exactly one point, denoted $\alpha $ in (\cite{bretechetorique}, Lemme 4.3), lies in the polar divisor of $Z(U(\bold
A);\s).$  This   is proved by showing that the function   
$G(\s)$ [ibid. (4.2)], satisfies $G(\bold 0) \neq 0.$ The proof of this property is actually 
indirect, and does not follow from the fact that $G$ is defined at $\bold 0.$ This is because $G$ has  both  positive and
negative coefficients in its series expansion at $\bold 0.$   In our notation, $G(\s) = {\cal H_\alpha} (\s).$ For us, the
fact that $G(\bold 0) \neq 0$   is a very   special  case of the {\it general property} 
(\ref{nonzero}) that applies to  {\it  any point in} $\partial \left( K \cup S^*(W) \right)^o\,.$ The proof of
(\ref{nonzero})   is both direct and independent of any need to desingularize $X(\bold A).$ This also enables us to work with
concrete examples.\hfill
$\square$

\medskip

To proceed, we
will   need to introduce some additional notations, and prove a preliminary result. First, we write   $W$ as a polynomial
by setting
$W(X_1,\dots,X_n)=1 + \sum_{\bold \nu \in S^*(W)} u(\bold \nu) X^{\bold \nu}.$  

In addition, we define $I = K \cup S^*(W).$ Since this is a finite set, we have (see Notations)   $\tilde I = I,$ so that
$\tilde I^o = I^o = \{\x \in \R_+^n : \la \x, \nu \ra \ge 1 \ \ \forall \nu \in I\}.$ We set  $\Gamma = \partial I^o\,,$ and
for any $\alpha \in \Gamma,$ define   $ K(I, \alpha) = \{\nug \in I: \la \alpha, \nug \ra = 1\}. $  
 

Finally, for all $\bold \nu \in K \cup S^*(W),$ we  define $c'(\bold \nu)$ as follows:
\be
\item $c'(\bold \nu)=c(\bold \nu)$ if $\bold \nu \in K\setminus S^*(W)$;
\item $c'(\bold \nu)=u(\bold \nu)$ if $\bold \nu \in S^*(W) \setminus K$;
\item $c'(\bold \nu)=c(\bold \nu)+u(\bold \nu)$ if $\bold \nu \in K\cap S^*(W)$.
\ee
The following lemma  plays an important role in
the proof of Theorem \ref{application2bis}:
\begin{Lemme}\label{precision}
For each $\alpha \in \Gamma,$ and   each  $\bold \nu \in K(I;\alpha),$\  $c'(\bold \nu)=1$.
\end{Lemme}
{\bf Proof:}\par
We start  with the presentation (\ref{hA}), and  choose  
$\eta < \frac{1}{4}\min_{\bold \nu \in I\setminus K(I;\alpha)} 
\left( \la \alpha , \bold \nu \ra -1 \right)$ \newline {\it if} $K(I;\alpha)\neq I.$ Otherwise, we choose  $\eta \in
(0,1/6)$.\\ We set 
${\cal F}= \{\eps \in (0,1)^{2n} : 1,\eps_1,\dots,\eps_{2n} {\mbox { are 
linearly independant over }}\Q\}$.\\
For each $\eps \in {\cal F}$ we define:
\be
\item 
$\alpha(\eps)=\left(\alpha_1(\eps),\dots,\alpha_n(\eps)\right),$ \ where 
$\alpha_i(\eps)=(1-\eps_i)\alpha_i +\eps_{n+i}$ for all
$i=1,\dots,n$;
\item 
$g_\eps(t)=h_{\bold A} (t^{\alpha_1(\eps)},\dots,t^{\alpha_n(\eps)})$ 
for all $t\in (0,1)$.
\ee
By using the  bound for $\eta,$ as above, and the fact that $\la \alpha (\eps)  , \bold \nu \ra 
=\la \alpha , \bold \nu\ra + O(|\eps|)$ as $|\eps| \rightarrow 0$ \ (since $I$ is finite), 
it is clear that one can choose $\eps \in {\cal F}$ with $|\eps|$ so small that the following property is satisfied:\par
\quad $\bold \nu \in K(I;\alpha)$ \ {\it implies } \ $\la \alpha (\eps)  , \bold \nu \ra <1+\eta$ \
{\it and }
\begin{equation}\label{devlimite}
g_\eps(t) = 1+\sum_{\bold \nu \in K(I;\alpha)} c'(\bold \nu) \, t^{\la \alpha (\eps)  , \bold \nu \ra}
  +O_\eps(t^{1+\eta}) \quad (t\rightarrow 0).
\end{equation}

We fix any such  $\eps$ in the following. 

On the other hand, it is also clear that there exist $N=N(\eta,\eps)$ 
such that 
\begin{equation}\label{devlimitebis}
g_\eps(t) = \sum_{\bold \nu \in T({\bold A})\atop |\bold \nu|\leq N } 
t^{\la \alpha (\eps)  , \bold \nu \ra} +O_\eps(t^{1+\eta}) \quad
(t\rightarrow 0).
\end{equation}
Since $\epsilon \in  {\cal F},$  it follows that if $\bold \nu \neq \bold \nu' \in \N_0^n,$ then 
$\la \alpha (\eps)  , \bold \nu \ra 
\neq \la \alpha (\eps)  , \bold \nu' \ra.$ In particular, this insures that for any $\bold \nu \in K(I,\alpha),$ the
coefficient of
$t^{\la \alpha (\epsilon),\bold \nu \ra}$ in (\ref{devlimite}) equals $c'(\bold \nu),$ and in  (\ref{devlimitebis})
equals $1.$ Since the two partial asymptotic expansions must be equal  up to terms of order $t^{1+\eta},$ this shows that
$c'(\bold \nu) = 1$ if 
$\bold \nu
\in K(I,\alpha).$ \hfill   
$\square$. 

Our first observation is as follows. 
\begin{Theoreme}\label{application2bis} For each point   $\alpha \in \Gamma,$ the meromorphic continuation
of $Z(U(\bold A);\s)$ is not analytic at $\alpha.$ 
\end{Theoreme}
{\bf Proof:}\par

We need to sharpen the   proof of  Theorem 3. Let $\alpha \in \Gamma$ be arbitrary and fixed.  Using   the same
argument as in the proof of Theorem 3, we first check  that 
$Z(F_{\bold A};\s)$ converges absolutely if 
$\sigma_i >\alpha_i$ for each $i,$  

We next introduce the product of linear forms ${\cal L}_\alpha(\s):= \prod_{\bold \nu \in K(I;\alpha)}\, \la\bold
\nu,\s\ra,$ and use Lemma
\ref{precision} to write it as follows:
$$ {\cal L}_\alpha(\s)  =  \prod_{\bold \nu \in K \cap K(I,\alpha)} \la\bold \nu,\s\ra^{c(\bold \nu)} \,  \cdot 
\prod_{\bold \nu
\in  S^*(W) \cap K(I,\alpha)} \la\bold \nu,\s\ra^{u(\bold \nu)}\,.$$
The function ${\cal H}_\alpha(\s):=  Z(F_{\bold A};\alpha + \s) \cdot {\cal L}_\alpha(\s)$ is evidently analytic in 
$V(0)=\{\s\in
\C^n : \,\la
\nu ,\sigma \ra >0\ 
\forall
\nu
\in I\}.$ We first show that it is  analytic in some larger domain $V(-\delta_1)$ for some positive
$\delta_1,$ by grouping each   factor  in ${\cal L}_\alpha(\s)$ with an appropriate factor of $Z(F_{\bold A};\alpha + \s)$
obtained from (\ref{zetafunc}). 

For the leftmost factor on the rightside of (\ref{zetafunc}), we have:  
\begin{eqnarray*}
&{}&\prod_{\bold \nu\in K} 
\zeta\left(\la \bold \nu,\alpha\ra+\la \bold \nu,\s\ra\right)^{c(\bold \nu)}\,\cdot \, \prod_{\bold \nu \in K \cap
K(I,\alpha)}
\la\bold \nu,\s\ra^{c(\bold \nu)} \\ &=&  \prod_{\bold \nu \in K\cap K(I;\alpha)} \big[\la \bold \nu, \s \ra \cdot 
\zeta\left(1+\la \bold \nu,\s\ra \right)\big]^{c(\bold \nu)}\,\cdot \,   
 \prod_{\bold \nu\in K\setminus K(I;\alpha)} 
\zeta\left(\la \bold \nu,\alpha\ra+\la \bold \nu,\s\ra\right)^{c(\bold \nu)}. 
\end{eqnarray*}
For $\delta_0$ chosen small enough, it is clear that each of the two products on the last line, one over $\bold \nu \in
K\cap K(I;\alpha),$ the other over $\bold \nu \in K - K(I,\alpha),$  is analytic in $V(-\delta_0).$

For the rightmost factor on the right side of (\ref{zetafunc}),   observe first that (\ref{zetafunc})  and the proof of
Lemma 1 imply that there exists\, 
$\delta
\in (0, 1)$ such that 
\begin{equation}\label{gdelta}
G_\delta (\s) := Z(W; \s) \cdot \prod_{\bold \nu \in S^*(W) \cap K(I,\alpha)} \zeta (\la \bold
\nu,
\s
\ra)^{-u(\bold
\nu)}\ \ \ {\mbox {is analytic in $V(W; 1-\delta).$}}\end{equation}
Thus, 

$Z(W;\alpha + \s) \prod_{\bold \nu \in S^*(W) \cap K(I,\alpha)} \la \bold \nu, \s
\ra^{u(\bold \nu)} = \prod_{\bold
\nu\in S^*(W)\cap K(I,\alpha)} \big[\la \bold \nu, \s \ra \zeta(1 + \la \bold \nu, \s \ra)\big]^{u(\bold \nu)}\, \cdot \,  
G_\delta (\alpha +
\s),$ 

and   $G_\delta (\alpha + \s)$ is analytic for $\s \in V(-\delta_0'),$ for some $\delta_0' > 0.$

We conclude that   ${\cal H}_\alpha(\s)$ can be written in $V(0)$ as follows:
\begin{eqnarray*}
{\cal H}_\alpha(\s) &=&  \prod_{\bold \nu\in K\cap K(I;\alpha)} \big[\la \bold \nu, \s \ra \cdot 
\zeta\left(1+\la \bold \nu,\s\ra\right) \big]^{c(\bold \nu)}\  \cdot  
 \prod_{\bold \nu\in K\setminus K(I;\alpha)} 
\zeta\left(\la \bold \nu,\alpha\ra+\la \bold \nu,\s\ra\right)^{c(\bold \nu)} \\
&{}&\quad \ \ \ \cdot  \prod_{\bold \nu\in S^*(W)\cap
K(I,\alpha)} \big[\la \bold \nu, \s \ra \cdot \zeta(1 + \la \bold \nu, \s \ra)\big]^{u(\bold \nu)}\,   \cdot G_\delta
(\alpha +
\s)\\ &=& \prod_{\bold \nu \in K\cap K(I;\alpha)} \big[\la \bold \nu, \s \ra \cdot 
\zeta\left(1+\la \bold \nu,\s\ra\right) \big]^{c(\bold \nu)} \  \cdot  \prod_{\bold \nu\in S^*(W)\cap
K(I,\alpha)} \big[\la \bold \nu, \s \ra \cdot \zeta\left( 1 + \la \bold \nu, \s \ra \right)\big]^{u(\bold \nu)}   \\
&{}&\quad \ \ \ \cdot   \prod_{\bold \nu\in K\setminus K(I;\alpha)} 
\zeta\left(\la \bold \nu,\alpha\ra+\la \bold \nu,\s\ra\right)^{c(\bold \nu)} \,  \cdot G_\delta (\alpha +
\s),\end{eqnarray*}   and we know that there exists $\delta_1' > 0$ such that the product of the two functions on the
last line is  analytic in $V(-\delta_1').$ 

Applying    Lemma \ref{precision} a second time now shows that for any $\s \in V(0):$
\begin{eqnarray}\label{hform} 
{\cal H}_\alpha(\s) &=& \prod_{\bold \nu \in K(I,\alpha)} \big[\la \bold \nu, \s \ra \cdot 
\zeta\left(1+\la \bold \nu,\s\ra\right) \big] \\
&{}& \quad \ \ \cdot  \prod_{\bold \nu\in K\setminus K(I;\alpha)} 
\zeta\left(\la \bold \nu,\alpha\ra+\la \bold \nu,\s\ra\right)^{c(\bold \nu)}\,   \cdot G_\delta (\alpha + \s).\nonumber\end{eqnarray}
We then deduce the existence of $\delta_1 > 0,$ such that the product over $\bold \nu \in
K(I,\alpha)$ in the first line of (\ref{hform}) is analytic in $V(-\delta_1).$ Since the product of functions on the second
line   is analytic if
$\delta_1$ is chosen sufficiently small,  we have  verified what we needed to show,  that is,
${\cal H}_\alpha(\s)$ is analytic in some neighborhood $V(-\delta_1)$ containing $\s = \bold 0.$ 

The second part of the argument is   an immediate consequence of the following essentiel property:
\begin{equation}\label{nonzero} {\cal H}_\alpha(\bold 0) \neq 0.\end{equation}  
To prove this, we start with (\ref{hform})
and rewrite the product by writing
 
$\quad 1 = \prod_{\bold \nu \in K\cap K(I,\alpha)} \zeta(1 + \la \bold \nu, \s
\ra)^{c(\bold \nu)}
\cdot \prod_{\bold \nu \in K\cap K(I,\alpha)} \zeta(1 + \la \bold \nu, \s \ra)^{-c(\bold \nu)}.$  

Multiplying together all the
terms with exponent
$-c(\bold \nu)$ with the factor 
 $\prod_{\bold \nu \in S^*(W) \cap K(I,\alpha)} \zeta(1 + \la \bold \nu, \s \ra)^{-u(\bold \nu)}$ in (\ref{gdelta}), evaluated
at $\alpha + \s,$  and applying  Lemma \ref{precision} again, gives a factor of ${\cal H_\alpha}(\s)$ that equals \ 
$\prod_{\bold \nu \in K(I,\alpha)} \zeta(1 + \la \bold \nu, \s \ra)^{-1}.$
Multiplying together all the terms with exponent $c(\bold \nu)$ with the product over $\bold \nu \in K - K(I, \alpha)$ 
in (\ref{hform}) gives a factor equal to \  
$\prod_{m\in K} \zeta (\la \bold \nu, \alpha + \s\ra)^{c(\bold \nu)}.$ 
Thus, we find a different expression for ${\cal H}_\alpha(\s)$ as a product of functions, each of which is analytic, at
least, in
$V(0):$
\begin{eqnarray}\label{h2} 
{\cal H}_\alpha(\s) &=&   \prod_{m\in K} \zeta (\la \bold \nu, \alpha + \s \ra)^{c(\bold \nu)}  \cdot  \prod_{\bold \nu \in K(I,\alpha)}
\zeta (1 + \la \bold \nu, \s \ra)^{-1}   \cdot Z(W,\alpha + \s)  \\
&{}& \ \cdot \prod_{\bold \nu \in K(I,\alpha)} \big[\la \bold \nu, \s \ra \cdot 
\zeta\left(1+\la \bold \nu,\s\ra\right) \big]\,. \nonumber 
\end{eqnarray}

Since there exists a neighborhood of $\s = \bold 0$ in which the  function $\prod_{\bold \nu \in K(I,\alpha)} \big[\la \bold \nu, \s \ra
\cdot 
\zeta\left(1+\la \bold \nu,\s\ra\right) \big]$ is both analytic and {\it never } $0,$    it follows that the product   in 
(\ref{h2}) is
actually analytic in a neighborhood of $\s = \bold 0.$ In such a neighborhood, we therefore have:
\begin{equation}\label{prolongementh}
{\cal H}_\alpha(\s) = \prod_p H(p;\s)  \cdot  \prod_{\bold \nu \in K(I;\alpha)} 
\big[\la \bold \nu, \s \ra \cdot 
\zeta\left(1+\la \bold \nu,\s\ra\right) \big]\,, 
\end{equation}
where 
$$H(p;\s)=  
\prod_{\bold \nu \in K} 
 (1-p^{-\la \bold \nu,\alpha\ra-\la \bold \nu,\s\ra})^{-c(\bold \nu)}\,  \cdot \prod_{\bold \nu \in K(I;\alpha)}
{ (1-p^{-1-\la \bold \nu,\s\ra})}\,  \cdot
 W\left(p^{-\alpha_1 -s_1},\dots,p^{\alpha_n-s_n}\right)\,.$$
The function $\s \to \prod_p H(p;\s) $ is analytic at $\s = \bold 0,$ but we still need to understand its value at this
point.  For
$r\in (0,1)$  we  define the open neighborhood \
${\cal B}(r)=V(0)\cup
\{\s
\in
\C^n
\mid |s_i|<r\}$
\ of
$\bold 0,$ and
  write out   $ H(p;\s)\big|_{{\cal B}(r)}.$ For our purposes, it suffices to observe the   existence of $u > 1$ such
that  the following holds, to which we apply    Lemma \ref{precision} for the last
equation:
\begin{eqnarray*}
H(p;\s)&=&\bigg(
1+\sum_{\bold \nu \in K\cap K(I;\alpha)} 
\frac{c(\bold \nu)}{p^{1+\la \bold \nu,\s\ra}}+ O(p^{-u+r})\bigg) \bigg(
1-\sum_{\bold \nu \in K(I;\alpha)} 
\frac{1}{p^{1+\la \bold \nu,\s\ra}}+ O(p^{-u+r})\bigg)\\
& &
\cdot  \bigg(
1+\sum_{\bold \nu\in S^*(W)\cap K(I;\alpha)} 
\frac{u(\bold \nu)}{p^{1+\la \bold \nu,\s\ra}}+ O(p^{-u+r})\bigg) 
\\
&=&
1-\sum_{\bold \nu\in K(I;\alpha)} 
\frac{1-c(\bold \nu)-u(\bold \nu)}{p^{1+\la \bold \nu,\s\ra}}+ O(p^{-u + r })  = 1-\sum_{\bold \nu\in K(I;\alpha)} 
\frac{1-c'(\bold \nu)}{p^{1+\la \bold \nu,\s\ra}}+O(p^{-u+r })\\
&=&1+{\cal O}(p^{-u+r}) \quad {\mbox {uniformly in }} \s \in {\cal B}(r).
\end{eqnarray*}
Thus, by choosing $r$ so small that $-u + r < -1$ for all $\s \in {\cal B}(r),$ we conclude that   $\s \mapsto \prod_p H(p;\s)$
also converges absolutely in ${\cal B}(r).$    We can therefore evaluate both sides of
(\ref{prolongementh}) at $\s = 0.$ In this way, we find the following Euler product expansion that {\it converges} to
${\cal H}_\alpha(\bold 0):$
\begin{equation}\label{hvalue}
{\cal H}_\alpha(\bold 0) = \prod_{p} \bigg( (1-p^{-1})^{\#K(I,\alpha)}\, \cdot W(p^{-\alpha_1},\dots, p^{-\alpha_n})\,
\cdot
\prod_{\bold \nu
\in K} \big(1-  p^{-\la \bold \nu,\alpha \ra}\big)^{-c(\bold \nu)} \bigg)\,.\end{equation}
The distinct advantage of (\ref{hvalue}) is that it easily is seen to imply that ${\cal H}_\alpha (\bold 0) > 0.$ Indeed,
we know that $$W(p^{-\alpha_1},\dots, p^{-\alpha_n})\, \cdot
\prod_{\bold \nu
\in K} \big(1-  p^{-\la \bold \nu,\alpha \ra}\big)^{-c(\bold \nu)} = h_{\bold A}(p^{-\alpha_1},\dots, p^{-\alpha_n}) > 0 \ \ {\mbox {for
each
$p.$}}$$
Thus,    each factor of the Euler product in (\ref{hvalue})  is {\it positive}. This implies ${\cal H}_\alpha (\bold
0)$ is also positive. As a result, the equation that gives the meromorphic continuation of $Z(U(\bold A);\s)$ in a
neighborhood of $\alpha,$   
$$Z(U(\bold A);\alpha + \s) = \frac {{\cal H}_\alpha(\s)}{{\cal L}_\alpha(\s)}\,,$$
now implies that the right side {\it cannot}  be analytic at $\s = \bold 0.$ This completes the proof of Theorem
\ref{application2bis}.\hfill$\square$  

\medskip

To state the second main result, we first need some notions from \cite{lichtinduke1}. Since $s_1\cdots s_n$ divides
$Z(F_{\bold A};\s)$ in (\ref{perron}), we will work with the extended polyhedron, by setting $X = I \cup \{\bold
e_j\}_1^n $ and $\hat \Gamma = \partial X^o.$ By definition, a vertex of $\hat \Gamma$ is the intersection of $n$
linearly independent support planes to $\hat \Gamma.$ Set ${\cal V}$ to denote the set of vertices of $\hat \Gamma.$ For
each
$\bold \alpha
\in {\cal V},$ there is an $n$ dimensional  closed cone ${\cal C}(\bold \alpha)$ of direction vectors in $(0,\infty)^n$ 
defined by the property:
$$\bold \gamma \in {\cal C}(\bold \alpha) \quad {\mbox { iff }} \quad \big\{\sigma \in \R^n: \la \bold \gamma, \sigma \ra
=
\la
\bold
\gamma,
\bold
\alpha
\ra\big\}
\ {\mbox {is a support plane of
$\hat \Gamma$}}.$$
Any vector in the interior of ${\cal C}(\bold \alpha),$ for some vertex $\bold \alpha,$  is called a {\it generic}
(direction) vector. It is clear that the set of generic vectors is an open dense subset of $(0,\infty)^n.$ 
To each vertex $\bold \alpha$ of $\hat \Gamma,$ there exists the subset  $\hat K(I, \bold
\alpha) =
\{\bold \nu_1,\dots, \bold \nu_m\} \subset X,$ \ $m = m (\bold \alpha)\ge n,$ such that ${\it {rank}} \{\bold \nu_i\} = n,$\, 
and the polar locus of
$\hat Z(U(\bold A);\s) := Z(F_{\bold A};\s)/s_1\cdots s_n$ through $\bold
\alpha$ is the union of affine  planes \ $\bigcup_{ \bold \nu_i \in \hat K(I, \bold
\alpha) } \,\big\{\la \bold \nu_i, \s \ra = \la \bold \nu_i, \bold \alpha \ra 
\big\}\,.$ 
\begin{Theoreme}\label {coefficient} Let $\bold \alpha $ be a vertex of \ $\hat \Gamma,$ and $\gamma $   a
generic vector in ${\cal C}(\bold \alpha).$ Then there exists a nonzero polynomial $Q_\gamma (u)$ of degree $m(\bold
\alpha) - n,$ and some
$\theta > 0,$  such that:
$$N_\gamma (U(\bold A), t) := C(\bold A) \cdot \sum_{(m_1,\dots, m_n)\in N^n \atop m_i \le t^{\gamma_i}\ \forall i}
F_{\bold A}(m_1,\dots, m_n) =   t^{\la \gamma, \bold \alpha \ra}\, Q_\gamma (\log t) + O\left(t^{\la \gamma, \bold \alpha
\ra -
\theta}\right)\ \ {\mbox {as }} \ t \to \infty.$$  
\end{Theoreme}
{\bf Proof:} The appropriate analog of ${\cal H}_\alpha (\s)$ when $\hat Z(U(\bold A);\s)$ replaces
$Z(U(\bold A);\s)$ in the proof of Theorem \ref{application2bis} is, in the preceding notation,  given by  $\hat {\cal
H}_\alpha (\s): =
\hat Z(U(\bold A);\alpha +
\s)
\cdot
\prod_{\bold \nu_i \in \hat K(I, \bold \alpha)} \la \bold \nu_i, \s \ra.$ The proof of 
the fundamental fact (\ref{nonzero})  extends straightforwardly to show $\hat {\cal H}_\alpha (\bold 0) \neq 0.$  This now 
allows us to  apply  Th\'eor\`eme 2 part iv of
\cite{bretechecompo}   since the constant
$C_0$ in the notation of [ibid., 1.10] equals, in our notation,   $\hat {\cal H}_{\bold \alpha}(\bold 0).$  
It should also be noted that the proof in [ibid.]   gives an explicit  expression for $Q_\gamma (\log t)$ as a certain
volume integral. This however is not needed for   purposes of this article. \hfill $\square$. 

{\bf Remark:} It would be interesting to know if the argument in \cite{bretechetorique} could extend to prove Theorems
\ref{application2bis}, \ref{coefficient}, but we do not see how to prove the crucial nonvanishing result (\ref{nonzero})
using  the methods in [ibid.] that exploit  Mo\"ebius inversion.  

\subsection{How often is the product of $n$ integers an 
$n^{th}$ power?}

A natural problem in multiplicative number theory is to describe the asymptotic density  of $n-$fold products of positive
integers that also equal the $n^{th}$ power of an integer. When $n = 3,$ several authors have given a precise
asymptotic for the density. Their starting point was an observation of Batyrev-Tschinkel
(\cite{salberger},11.50) who noted that the problem is equivalent to finding the asymptotic of the exponential height density
function on a certain singular cubic toric variety.  This interpretation naturally extends to any   $n \ge 3.$ However, until
now, no extension of these results  to arbitrary $n$ seems to have been published  in the literature. The purpose of this
subsection is to solve the problem for arbitrary $n \ge 3$ by applying a variant of the  methods of \S 2.3. A point that must
be addressed is the fact that Theorem
\ref{coefficient} only applies to generic directions. However, for this problem the direction $(1,\dots, 1)$ is of special
interest, and it is {\it not} a priori clear that this is   generic.

In the following discussion, we use the notations from \S 2.2, 2.3. In particular,   
$\bold A_n = (1,\dots, 1, -n)$  is the appropriate $1 \times (n+1)$ integral matrix whose row sums to $0.$  Note that
$U(\bold A_n)$ is now defined to be:
$$U(\bold A_n) = \{\x = (x_1:\cdots:x_{n+1}) \in \P^n(\Q) :  x_1\cdots x_n = x_{n+1}^n \ \ {\mbox{ and }}\  x_1\cdots x_n
\neq 0\}.$$ The density zeta function of interest is
$$ Z(U(\bold A_n);\s) := \sum_{\x\in U(\bold A_n)} H(\x;\s) \quad {\mbox {where}} \ \    \s=(s_1,\dots,s_{n+1}), $$ 
and $H(\x;\s)$ is defined as in \S 2.3, using the  unique integral vector representative of a point $\x$ with 
components whose $gcd$ equals $1.$ Setting $\r = (r_1,\dots, r_n),$ we also define
\begin{eqnarray*}
D_n &=& \bigg\{\r \in \{0,\dots,n-1\}^n :\   
 \frac{ |\r|}{n} \in \N \bigg\},\  \ {\mbox { where }} \ |\r| = r_1 + \cdots + r_n,\\ 
J_n &=& \bigg\{\r+\e_{n+1} : \r \in 
\{0,\dots,n\}^n  {\mbox { and }} |\r|=n\bigg\} \setminus \{(1,\dots,1)\},\\
\ell (\r) &=& \r + \frac{|\r|}{n}\e_{n+1} = \big(r_1,\dots, r_n, \frac{|\r|}n\big) \quad {\mbox { for any } }\r \in D_n\,,
\end{eqnarray*} and for every $\delta \in \R$,  
$$V(\delta):= 
\{\s \in \C^{n+1}: \  
 \la \ell (\r), \sigma \ra  >\delta  \ \ \forall \r \in D_n\}.$$
\begin{Theoreme}\label{application3} For any $n \ge 3$ the following three assertions are satisfied.
 \be
\item 
$\s \mapsto Z\left(U(\bold A_n);\s\right)$ converges absolutely
in $V(1)$ and satisfies :
$$Z\left(U(\bold A_n);\s\right)=\frac{\prod_{i=1}^n \zeta
  (ns_i+s_{n+1})}{\zeta(s_1+\dots+s_{n+1})}\,\cdot   
\prod_{p} \bigg(\sum_{\r \in D_n} 
\frac{1}{p^{\la \ell(\r), \s \ra}} \bigg);$$
\item 
$\s \mapsto Z\left(U(\bold A_n);\s\right)$ can be
meromorphically continued to $V(0)$ and $\partial V(0)$ is the natural boundary of $Z(U(\bold A_n);\s);$
\item 
there exists $\theta >0$
such that:
$$ N_\infty(U(\bold A_n);t) =  t  Q_n(\log t ) +O(t^{1 -\theta}) {\mbox { as }} t\longrightarrow \infty,$$
where $Q_n$ is a non-vanishing polynomial of degree  $d_n = {2n-1 \choose n} -n-1$ satisfying
$$ Q_n(\log t  ) =C_0(n)~t^{-1} Vol(A_n(t))  +
O({\log^{d_n - 1} (t)})\ \  {\mbox { as }}\  t\longrightarrow \infty,$$ 
\ee
 $A_n(t)$ is defined with the help of the vector $\beta: = \left(1,\dots,1,1+\frac{1}{d_n +
1}\right)$ to equal  
\begin{eqnarray*}
A_n(t)&=& \bigg\{\x =(x_\nu)_{\nu \in J_n} \in
{[1,+\infty[}^{d_n + n} \, :  
\prod_{\nu \in J_n} x_\nu^{\nu_j} \leq t^{\beta_j} \quad \forall j=1,\dots,n+1 \bigg\},\\
{\mbox{and}} \quad \ \  C_0(n) &=& 2^{n-1} \cdot \prod_{p}\bigg(
{(1 - p^{-1})}^{d_n + 1}\,\cdot  \sum_{\r \in D_n } p^{-\frac{|\r|}{n}}\bigg) >0 \end{eqnarray*}  

\end{Theoreme}

{\bf Proof:} Defining  
$$  T (\bold A_n)=\{\alpha \in \N_0^{n+1} : 
\alpha_1+\dots+\alpha_n=n\alpha_{n+1} {\mbox { and }} \alpha_1 \dots
\alpha_{n+1}=0 \},$$  
we first need to   construct an explicit presentation of   
$$h_{\bold A_n} (X) = \sum_{\alpha \in T (\bold A_n) }
 X_1^{\alpha_1}\dots X_{n+1}^{\alpha_{n+1}}.$$
To do so, we observe that for every $X \in P(1)$:
\begin{eqnarray*}
h_{\bold A_n} (X)
&=&
\sum_{\alpha_1+\dots+\alpha_n=n\alpha_{n+1} \atop  
\alpha_1 \dots \alpha_{n+1}=0} X^\alpha 
= (1-X_1\dots X_{n+1}) \cdot \sum_{\alpha_1+\dots+\alpha_n=n\alpha_{n+1}} X^\alpha \\
&=& (1-X_1\dots X_{n+1}) \cdot \sum_{n \big| \alpha_1+\dots+\alpha_n}
X_1^{\alpha_1}\dots X_n^{\alpha_n} 
X_{n+1}^{\frac{\alpha_1+\dots+\alpha_n}{n}}\\
&=& (1-X_1\dots X_{n+1}) \cdot \sum_{\r \in D_n}
X_1^{r_1}\dots X_n^{r_n} 
X_{n+1}^{ |\r|/n}\, \cdot 
\sum_{\alpha \in \N_0^n}
X_1^{n\alpha_1}\dots X_n^{n\alpha_n} 
X_{n+1}^{|\alpha|}\\
&=& \bigg(\prod_{i=1}^n {(1-X_i^n
    X_{n+1})}^{-1}\bigg) \cdot W_n(X_1,\dots,X_{n+1}).
\end{eqnarray*}
We conclude that $(K,\la c(\bold \nu)\ra_{\bold \nu \in K},W_n)$ is a presentation of  
$h_{\bold A_n}(X)$ where:
\begin{eqnarray*}
 W_n(X_1,\dots,X_{n+1})&=& (1-X_1\dots X_{n+1}) \cdot \sum_{\r \in D_n}
X_1^{r_1}\dots X_n^{r_n} X_{n+1}^{ |\r| / n }\\ 
K &=& \{ n\e_i +e_{n+1} : i=1,\dots,n\} \\
 c(\bold \nu) &=& 1 \quad \forall \bold \nu \in K.\end{eqnarray*}
  
Assertion  1 and the first part of Assertion 2 of the Theorem now follow immediately from 
Theorem \ref{application2}.\par
To prove that $\partial V(0)$ is the natural boundary of $Z(U(\bold A_n);\s),$  it suffices to
show that the polynomial $W_n$ is not cyclotomic when $n \ge 3.$ We show this by contradiction. 

Thus, suppose that $W_n$ is
cyclotomic.  It is then clear that the polynomial 
$$W_n^*(X_1,\dots,X_{n+1}):=\sum_{\r \in D_n}
X_1^{r_1}\dots X_n^{r_n} X_{n+1}^{ |\r| / n } $$ is also cyclotomic. From this it follows  that
the polynomial in  one variable $R(t):=W_n^*(t,t,0,\dots,0,1)=1+(n-1)t^n$ is
cyclotomic. But this is impossible  since $R(t)$ has roots of modulus
 different from $1$. This completes the proof of Assertion 2.

\paragraph{Proof of Assertion 3}:\par
Setting  $I=K\cup S^*(W_n),$ elementary  computations  show the following properties:
\be 
\item
$W_n(X)=\sum_{ {}\atop {\r \in D_n  \atop \r \neq (1,\dots,1)}}
X_1^{r_1}\dots X_n^{r_n} X_{n+1}^{ |\r| / n }
-\sum_{{} \atop {\r \in D_n \atop \r \neq (0,\dots,0)}}
X_1^{r_1+1}\dots X_n^{r_n+1} X_{n+1}^{\frac{|\r|}{ n }+1}$;
\item $\iota (I)=1$, $\alpha^*=(\frac{1}{n},\dots,\frac{1}{n},0)=
\frac{1}{n}(\e_1+\dots+\e_n) \in R(I)$, $J(\alpha^*)=\{\e_{n+1}\}$
and \newline 
$K(I;\alpha^*)=J_n;$
\item $Rank\left(K(I;\alpha^*)\cup J(\alpha^*)\right)=n+1$ \  and \ 
$\# K(I;\alpha^*)= {2n-1 \choose n} -1 = d_n + n\,;$ 
\item the constant $C(\bold A_n),$ (see (\ref{zetarel})), equals $2^{n-1}.$
\ee
The vector $(1,\dots, 1)$ may not satisfy all the criteria needed to apply Th\'eor\`eme 2 part iv of \cite{bretechecompo}.
The idea is to find an equivalent vector as follows. Setting
$\beta =\left(1,\dots,1, 1+\frac{1}{d_n + 1}\right),$  it is easy to see that $\forall \x=(x_1,\dots,x_{n+1}) \in
\N^{n+1}$ satisfying 
$(x_1:\dots:x_{n+1})\in U(\bold A_n)$ and $\gcd(x_1,\dots,x_{n+1})=1$ 
we have 
$$ \max_i \,x_i \leq t \Longleftrightarrow 
 x_j \leq t^{\beta_j} \ \quad \forall j = 1,\dots, n+1,\ \ \ \forall t \ge 1.$$ 
To finish the proof, it suffices  to verify the criterion of [ibid.] that  there exists  
$\{\gamma_\nu\}_{\nu \in J_n \cup\{\e_{n+1}\}} \subset  (0,\infty)$
such that $\beta = \sum_{\nu \in R(I;\alpha^*)\cup J(\alpha^*)} \gamma_\nu
~\nu$.\par
We define: 
\be
\item $t(n)=\#\left\{\r \in \{0,\dots,n-1\}^n \, : \  
|\r| = n \right\}= d_n + 1$;
\item $\gamma_\nu ={t(n)}^{-1}$ \ \, $\forall \nu \in  \left(J_n \cup \{\e_{n+1}\}\right) \setminus \{n\e_i
  +\e_{n+1}\}_{i=1}^n\,,$ \ 
and \newline $\gamma_{n\e_i +\e_{n+1}} = 1 / n t(n)$ \ \ \ $\forall i=1,\dots,n. $  
\ee
We  first notice that   the value of $\sum_{{} \atop {\r \in \{0,\dots,n-1\}^n\, \atop |\r|=n}} r_j$ is independent of $j.$
Thus, for each $j=1,\dots,n$:
\begin{eqnarray*}
\sum_{\r \in \{0,\dots,n-1\}^n \atop |\r|=n} r_j
&=&
\frac{1}{n} \sum_{i=1}^n \ \sum_{\r \in \{0,\dots,n-1\}^n 
\atop |\r|=n} r_i
=\frac{1}{n}\ \sum_{\r \in \{0,\dots,n-1\}^n \atop |\r|=n} |\r|\\
&=&\#\left\{\r \in \{0,\dots,n-1\}^n\, :\  |\r| = n\right\}=t(n).
\end{eqnarray*}
A straightforward computation then shows:
\begin{eqnarray*}
\sum_{\nu \in R(I;\alpha^*)\cup J(\alpha^*)} \gamma_\nu \nu 
&=& \sum_{\nu \in J_n \cup \{\e_{n+1}\}} \gamma_\nu \nu \\
&=& (1+{t(n)}^{-1})\e_{n+1} + \sum_{i=1}^n {t(n)}^{-1}
\bigg(\sum_{\r \in \{0,\dots,n-1\}^n \atop |\r|=n} r_i \bigg)\e_i\\
&=& (1+{t(n)}^{-1}) \e_{n+1} + \sum_{i=1}^n  \e_i = \beta.
\end{eqnarray*}
This  completes the proof of
Theorem \ref{application3}.  

\section{Some applications in group theory}

The first two subsections give some simple applications of Theorems 1, 2 to two problems in the study
of a group zeta function that were first addressed by duSautoy and Gr\"unewald in  \cite{dSG00}, \cite{dSG02}. The third
section indicates an  additional application to a somewhat different subgroup counting problem that originates within the
theory of   finite abelian groups.  Recall that to a group $G,$ the
group zeta function is defined as follows:
$$\zeta_G (s) = \sum_{H \le G} |G : H|^{-s}.$$
\subsection {The largest pole of a cone integral}


The article \cite{dSG00} studied the group zeta function for a finitely generated nilpotent group $G.$ Its main result   was the
following.
\begin{Theoreme}\label{Theorem 3.1}  There exist a rational number $\alpha (G)$ and $\delta > 0$ such that $\zeta_G$ has
its largest real pole at $\alpha(G),$ and is meromorphic  in the halfplane $\sigma   > \alpha(G) - \delta.$
\end{Theoreme}
 
The proof given in [ibid] has two parts. First, $\zeta_G(s)$ is expressed in terms of an Euler product of ``normalized cone
integrals", whose $p^{th}$ factor (for a generic $p$) is analyzed by using ideas of Denef. The second part then uses methods
from the analysis of Artin L-functions to show the existence of a meromorphic continuation to the left of a first
(rational) pole.

The purpose of this subsection  is to show how Theorem 1 gives an alternative and more elementary proof of the second 
part of
the proof of Theorem \ref{Theorem 3.1}. Theorem \ref{Theorem 7}  is the essential part of this simpler argument. Thus, we
show that the  fundamental   result of \cite{dSG00},   the rationality of the abscissa of convergence of the group zeta function
for finitely generated nipotent groups,  can be proved by combining the   work of Denef
with the   methods of this paper. In addition, our method also proves a meromorphic
continuation of the group zeta function into a halfplane that contains its first real pole.

For the reader's convenience, we adopt the notation used in [ibid].  
 
We start with the representation of $\zeta_G(s)$ as an Euler product of normalized cone integrals: 
$$\zeta_G(s) = \prod_p a_{p,0}^{-1}\, Z_D (s - d, p),$$ 
where $d = $ Hirsch length of $G,$ and 
 $D = \{f_0, g_0, f_1, g_1,\dots, f_l, g_l\}$ $ \subset \Q [x_1,\dots, x_m]$ specifies  the cone integral data.   
 
Using work of Denef, it follows that for all sufficiently large $p,$ each $Z_D (s, p)$ can be expressed in a  purely
geometric-analytic fashion by using numerical data, produced by an embedded resolution of singularities $h: Y \to \Q^m$ for
the polynomial $F =
\prod_i f_i g_i$ (that is, of the reduced scheme $X: = {\it spec} (\Q[x_1,\dots, x_m]/(F)$), as follows:\par
$$Z_D (s, p) = \sum_{I \subset T} c_{p,I}\, P_I(p^{-s}, p),$$
where  $T$ denotes an index set for the irreducible $\Q$ components of $h^{-1}(X),$ and for each nonempty $I \subset T:$
\begin{eqnarray*}
c_{p,I} &:=& {\it card} \{a \in \overline {Y(\F_p)}: a \in \overline {E_i} \ {\it iff }\  i \in I\};\\
 \overline Y &:=&   \mbox {reduction of}\   Y  \ {\it { mod }} \  p;\\
P_I(p^{-s}, p) &:=& \frac{(p-1)^{|I|}}{p^m} \cdot \prod_{j\in I} \frac{ p^{-(A_j s + B_j)}}{1 - p^{-(A_j s +
B_j)}}.\end{eqnarray*}  

A characterization of the nonnegative integers $A_j, B_j$ is given in [ibid]. For our purposes here, it suffices to know that
to each divisor $E_j, j \in T,$ there corresponds a pair $(A_j, B_j)$ of nonnegative integers, at least one of which is
positive. It is also to be understood that we restrict attention to those $I$ for which $c_{p, I} > 0.$ The presence of this
factor indicates that $\zeta_G(s)$ is not ``uniform" (see Introduction). Although the expression of the product of the factors
$P_I$ given in [ibid.] is a priori more intricate, this is not really  needed to prove Theorems \ref{Theorem 3.1}  or
\ref{Theorem 7}.

The constant term $a_{p,0}$ of $Z_D (s,p)$ is independent of $p^{-s}$ and expressed as follows: 
$$a_{p,0} = p^{-m} \cdot \sum_{\{I: A_j = 0 \ \forall j \in I\}} (p-1)^{|I|}\,\cdot c_{p,I} \cdot \prod_{j\in I}
\frac{p^{-B_j}}{1 - p^{-B_j}}\,.$$

An important first observation is that $a_{p,0} > 0.$ Given this, it is then necessary to   bound each $c_{p,I}/a_{p,0}.$
In general, there is significant fluctuation in
$c_{p,I}$ as a function of $p$ (whence the ``nonuniform" nature of the Euler product). So, one cannot, as yet, hope to do
better than the following, which is a modest improvement over that proved in [ibid].

Set $d_I = m - |I|.$
\begin{Lemme}\label{Lemma 3.2} For each $I,$ there exists $\delta_I \ge 1/2$ such that for all $p$ sufficiently large
$$\frac{c_{p,I}}{a_{p,0}} = p^{d_I} \big( 1 + O(p^{-1 - \delta_I}) \big).$$\end{Lemme}
{\bf Proof:} duSautoy-Gr\"unewald show that there exists $T$ and $P_0$ such that $p \ge P_0$ implies:\par
$a_{p,0}^{-1} \le (1 - T p^{-1})^{-1}.$

It is also clear that $(1 - T p^{-1})^{-1} = 1 + O(p^{-1})$ if $P_0$ is sufficiently large. 

Next, one uses an argument of Katz in the Appendix of \cite{Hoo}, to justify the existence of an integer $v_I \in [1, d_I - 1]$
so that:\par
$c_{p,I} = p^{d_I} \big(1 + O(p^{\frac{v_I - d_I}2})\big).$

Setting $\delta_I = \frac{v_I - d_I}2$ finishes the proof.\hfill $\square$

Some notations will now be useful. For each $I$ and $k \in T,$  set:
\begin{eqnarray*}
A_I &=& \sum_{j\in I} A_j, \quad B_I = \sum_{j\in I} B_j, \quad l_I (s) = A_I s + B_I;\\
l_k (s) &=&  A_k s + B_k \ \ \ {\mbox { for each }} k \in T;\\
\alpha_k &=& \frac {1 - B_k}{A_k} \quad {\mbox { if }} A_k > 0, \quad {\mbox { and }}\quad  \alpha_k = -\infty \ \ {\mbox {if
not;}}\\
\alpha_0 &=& \max_k \{\alpha_k\}, \qquad m_0(I) = \# \{j \in I: \alpha_j = \alpha_0\}, \quad m_0 (D) = \# \{k \in T: \alpha_k
= \alpha_0\}. \end{eqnarray*}

The following lemma is elementary, and is left to the reader to verify as a straightforward exercise. We implicitly assume
that $p \ge P_0$ for a suitably chosen $P_0.$ 
\begin{Lemme} \label{Lemma 3.3} For $P_0$ sufficiently large, and for each $\theta \in (0, 1),$ there exists $q =
q(\theta)$ such that $\sigma > \alpha_0 - \frac 1q$ imply  the following   properties for  any $I \subset T:$

i) $\prod_{j\in I} \frac 1{1 - p^{-l_j (s)}} = 1 + O\big( p^{-(1 - \theta)} \big);$

ii) if $|I| > 1,$ or $|I| = 1$ and $m_0 (I) = 0,$ then \ $p^{-l_I (\sigma)} < p^{-(1 + \theta)};$

iii) setting ${\cal I} = \{I:  |I| > 1, \ {\mbox {or}} \ \     |I| = 1 \ {\mbox { and }}\  m_0 (I) = 0\},$
$$ \sum_{ I \in {\cal  I}}\, p^{-l_I(s)} \,\cdot \prod_{j\in I} \frac 1{1 - p^{-l_j(s)}} = O \big( p^{-(1 + \theta)} \big);$$

iv) $$\sum_{\{k: \alpha_k = \alpha_0\}} \,\frac{p^{-l_k(s)}}{1 - p^{-l_k(s)}} = \sum_{\{k: \alpha_k = \alpha_0\}}\, p^{-l_k(s)}
+ O\big( p^{-(1 + \theta)} \big).$$
\end{Lemme}

We deduce from this lemma the following equation that holds for each $\theta \in (0, 1)$ and $\sigma > \alpha_0 - \frac
1{q(\theta)}:$
\begin{equation}\label {14} \prod_{p \ge P_0} a_{p,0}^{-1}\, Z_D (s, p) = \prod_{p \ge P_0} \bigg(1 + \sum_{\{k: \alpha_k =
\alpha_0\}}
\, p^{-l_k(s)} + O\big( p^{-(1 +
\theta)}
\big) \bigg).\end{equation}

We then multiply  both sides of (\ref{14}) by $\prod_{\{k: \alpha_k = \alpha_0\}} \zeta (l_k(s))^{-1}\,,$ and apply the
reasoning in the proof of Theorem 1 to conclude as follows.
\begin{Theoreme}\label{Theorem 7} The Euler product  in (\ref{14})  has a pole of order $m_0(D)$ at $s = \alpha_0 \in \Q,$
and admits a meromorphic continuation into a halfplane $\sigma > \alpha_0 - \delta,$ for some $\delta > 0,$ with
$\alpha_0$ as its only pole.  Furthermore, there is at most polynomial growth as $|\tau| \to +\infty$ within this
halfplane.
\end{Theoreme}

This now suffices to complete the proof of Theorem \ref{Theorem 3.1} since it is clear that  each of the remaining
finitely many factors,
$Z_D(s, p),\ p < P_0,$ admits a meromorphic continuation to $\C$ whose poles have rational valued real parts. On the
other hand, it does not yet seem possible to prove that $\zeta_G$ itself has polynomial growth within this halfplane. The
reason for this is that its largest real pole, say $\alpha (G),$  {\it could} originate from one of the factors indexed by
some
$p < P_0$ (where there is {\it not} good reduction in Denef's sense). In this event,   there will be   infinitely many
poles  of a  term of the form $(1 - p^{-A (s - d) - B})^{-1},$ where $\alpha (G) = d - B/A.$ The real part of each
such pole (in $s$) evidently equals $\alpha (G).$ Such a function could not, evidently,  have moderate growth in any
unbounded vertical   strip containing $\alpha (G).$   

\subsection{The natural boundary of a uniform zeta function}


The behavior of $\zeta_G(s)$ is most easily understood when it is uniform, and equals an absolutely convergent Euler product in
some halfplane \ $\zeta_G(s) = \prod_p \zeta_{G,p}(s).$ Many examples have  been found where these properties are known to
occur (see  \cite{dSG02}, \cite{dusautoy}). Uniformity essentially means that there exists a
single polynomial (sometimes also allowed to be a rational function) $h = 1 + H(X_1, X_2) \in \Z [X_1, X_2], H(\bold 0) = 0,$
and cyclotomic function
$u(X_1, X_2)$ (i.e. a finite product of integral powers of cyclotomic polynomials)  such that $\zeta_{G,p} = u(p, p^{-s})
\cdot  h(p, p^{-s}).$ Thus, there exist finitely many integers $a_i, b_i, \epsilon_i$ such that $$\zeta_G(s) = \prod_{i=1}^M
\zeta(a_i s + b_i)^{\epsilon_i}\, Z(h;s), \qquad {\mbox {where }}\ \  Z(h;s) = \prod_p h(p, p^{-s}).$$  
The nontrivial behavior of $\zeta_G(s)$ is therefore found in that of $Z(h;s).$

For an integral polynomial $h$  as above, \cite{dusautoy} showed the existence of a meromorphic continuation
of
$Z(h;s)$ up to a certain vertical line $\sigma = \beta_0$ that serves as a conjectured natural boundary of $\zeta_G(s),$ {\it
unless}
$h$ itself is also a cyclotomic polynomial. The method of du Sautoy requires one {\it first} to express $h$ as an infinite
product of cyclotomic functions $h = \prod_{(n,m)\in \N^2} (1 - X_1^m X_2^n)^{c_{n,m}}\,,\ c_{n,m} \in \Z,$ {\it before
beginning} the extension of
$Z(h;s)$ outside a halfplane of absolute convergence. 

An explicit expression is also given for $\beta_0$ in [ibid.]. One first
writes
$H = \sum_{j=1}^M\, H_j (X_1) X_2^j$ and defines $\deg_{X_1}  H_j = n(j).$ Then 
$$\beta_0 = \max_j \big\{ n(j)/j \big\}.$$ 

It is clear that Corollary 1 also applies to prove  the existence of a meromorphic continuation of $Z(h;s).$  
It should be evident to the reader that our method is   simpler because it does not require an a priori factorization of
$h$ as an infinite (in general) product  of cyclotomics. One can, instead,   begin with the expression for $h$ as a
polynomial. 

Using our procedure, a presumed
natural boundary is then given by
$\partial V^\# (h;0).$ To find this set, it suffices to write $H = \sum_k h_k(X_2) X_1^k\,,\ $ and define \ $\ord_{X_2} h_k =
m(k).$ Then 
$$\partial V^\# (h;0) = \{\sigma = \beta_1\}, \qquad {\mbox {where}}\ \  \beta_1 = \max_k \big\{ k/m(k) \big\}.$$
A simple check also shows that the same value is   obtained if one applies Corollary 2 to the Euler product $Z^*(h;\s) =
\prod_p h(p^{-s_1}, p^{-s_2}),$ and estimates  the presumed natural boundary by the set $\partial \big( V_2^\# (h; 0) \cap
\{s_2 = -1\}
\big),$  where   $V_2^\#(h;0)$ is notation for  the set denoted   $V^\#(h;0)$ in the statement of Corollary
2. 

It is then useful to observe the following.
\begin{Lemme}\label {lemma 3.1}   $\beta_0 = \beta_1\,.$ \end{Lemme}
{\bf Proof:} The equality does not seem to be completely obvious since  the two expressions above for $H(X_1, X_2)$   are
not symmetric in $X_1, X_2.$  Set
$k_1$ to be the largest index such that
$\beta_1 = k_1 /m(k_1) .$ We then observe that $\deg_{X_1} H_{m(k_1)} = k_1,$ that is, one has $n(m(k_1)) = k_1.$  If the
equation did not hold, then $\deg_{X_1} H_{m(k_1)} > k_1,$ which
implies the existence of an integer $l>k_1$ such that 
$(l,m(k_1)) \in S^*(h).$ Since $l > k_1$ and $\frac l{m(k_1)} > \beta_1,$ it
follows that \ $\frac{l}{m(l)} \leq  \beta_1 < \frac{l}{m(k_1)}\,.$ Thus,
$m(l) > m(k_1).$ However, the fact that $(l,m(k_1)) \in S^*(h)$ and $(l,m(l)) \in S^*(h)$ requires that $m(l) \le m(k_1).$
So, such an $l$ cannot exist. We then conclude that $\beta_1 \le \beta_0.$ 

Conversely, let $l_0$ be the smallest index such that $\beta_0 = n(l_0)/l_0.$ We then show that $m(n(l_0)) = l_0.$ Indeed, if
this were not the case, then $m(n(l_0)) < l_0$ must occur since $(n(l_0),l_0) \in S^*(h)$ implies $m(n (l_0)) \le l_0.$
Thus, there exists an integer $k < l_0$ such that $(n(l_0),k) \in S^*(h).$ It then follows that 
$\frac{n(k)}{k} \leq \beta_0 < \frac{n(l_0)}{k},$ which implies $n(k) <
n(l_0).$  However, since  $(n(l_0),k) \in S^*(h)$ and 
$(n(k),k) \in S^*(h),$ we must have $n(l_0) \le n(k),$ which is not possible. Thus, $m(n(l_0)) = l_0,$ which implies $\beta_0 \le
\beta_1.$ \hfill $\square$

{\bf Remark:}\  It may be of interest to the nonspecialized reader to see a concrete example of the preceding work. Let
$G$ be a
$\Q$-algebraic group, and $\rho$ a $\Q$-rational representation
$\rho : G \to GL_n$.
The local zeta function of $G$ is defined as 
$$\zeta_{G, \rho, p}(s)=\int_{G_p^+} \mid  \det(\rho(g))\mid_p^{-s}d\mu$$
where $G_p^+=G(\Q_p)\cap M_n (\Z_p)$ , $\mid.\mid_p$ denotes the p-adic valuation, 
and $d \mu$ is the normalised Haar measure on $G(\Z_p)$.
 The global zeta function is  the Euler product
$\zeta_{G, \rho}(s)=\prod_p \zeta_{G, \rho, p}(s).$

Consider the symplectic group $G = GSp_6$ of $6\times 6$  matrices $A\in GL_6(\Q)$ which satisfy the
equations $AJA^t=\lambda J$ where $A^t$ is the transpose of $A$ , $J$ the standard symplectic matrix
and $\lambda \in \Q^*$ . Let $\rho$ be the natural representation. The zeta function of $G$ 
\cite{igusa} has been found to equal
$$ \zeta_{G}(s/3)=\zeta(s)\zeta(s-3)\zeta(s-5)\zeta(s-6)\prod_{p} h(p,p^{-s}),$$
where $h(X,Y)= 1+(X+X^2+X^3+X^4)Y + X^5Y^2 $.
Corollary 1 shows that the presumed natural boundary is $\sigma = 4.$ 
In \cite{dSG02} it was proved   that this line $\Re (s)=4$ {\it is} the  natural boundary. 

\subsection{A refinement of the zeta function for abelian groups}

Let $G$ be a finite abelian group. We can write it as the direct sum 

$$G\simeq \Z/n_1\Z \oplus \Z/n_2\Z \oplus \cdots \oplus \Z/ n_r\Z $$

where $n_i\mid n_{i+1}$ . For a subgroup $H$ of $G$ and $a \in \N,$  
we consider the subgroup counting function  
$\tau_a(G)=\sum_{H<G}\, (\# H)^a\,,$  and define the zeta function associated to $\tau_a$ and $G$ with $r$ summands:
$$ Z^{(r)}(\tau_a, \s)=\sum_{\n \in \N^r \atop G\simeq \bigoplus_{j=1}^r \Z/n_j\Z }
\frac{\tau_a(G)}{n_1^{s_1}\cdots n_r^{s_r}}.$$ These zeta functions were studied in one variable  by specializing the values
of the $s_i$ to be 
$s_{r-k }=(k+1)s$ for $ 0\le k< r$. In particular it was proved in \cite {bhra2}
that
\begin{Theoreme}\label{ze}
  $ Z^{(r)}(\tau_a,s)$
is a rational function.
\end{Theoreme}
The explicit evaluation in the two variable case occurs in  \cite {bhra1}, i.e.
$$ Z^{(2)}(\tau_a,\s)=\zeta(s_1)\zeta(s_1-2a)\zeta(s_2)\zeta(s_2-2a)
\zeta(s_1-a-1)\prod_p\left(1+p^{a-s_1}-(p^a+1)p^{a-s_1-s_2}\right).$$
Applying Theorem 1 when $a=0,$ we obtain a meromorphic  continuation into the domain
$W(\l, 0)=\{\s : \sigma_1>0,   \sigma_1+\sigma_2)>0\},$ where $\l (\s) = (s_1, s_1 + s_2).$

Using an iteration of standard one variable Tauberian methods, we can then calculate an average order of the number of
subgroups as follows:
\begin{Theoreme} \label{th9}  $$\sum_{n_1\le x_1,n_2\le x_2 \atop G\simeq \Z/n_1\Z \oplus \Z/n_2\Z }\tau_0(G)=B
x_1^2x_2\log x_2+O(x_1x_2).$$ \end{Theoreme}

We compare this with the corresponding result (see \cite {bhwu}) in one variable in which the sizes of the constituent summands
of $G$ are not taken into account:
$$\sum_{ \# G \le x}\tau_0(G)=A_1x\log^2 x +A_2 x \log x +A_3 x +\Delta (x),$$
where $A_i$ are constants and $\Delta (x)\ll x^{5/8}\log ^4 x$.

Since the number of non-isomorphic abelian groups of order at most $x$ 
is aymptotically $x$, we see that the average number of subgroups in this case
is $\log^2 x.$ However,  when we also  consider  the size of each
direct summand, the average becomes   $x_1\log x_2$. Such greater   precision in the asymptotic behavior of  
  counting functions should be expected whenever the analytic study of multivariate zeta functions can be  combined with 
several variable Tauberian theorems such as those in \cite{bretechecompo}, \cite {lichtinduke1}, \cite {lichtinduke2},
\cite{lichtincomp}.

\bibliographystyle{plain}
\bibliography{bibliaeuler}

\begin{thebibliography}{10}

\bibitem{batyrev}
V.V. Batyrev and Yu. Tschinkel.
\newblock Manin's conjecture for {T}oric variety.
\newblock {\em Journal of Algebraic Geometry}, 7(1):15--53, 1998.

\bibitem{bhra1}
G.~Bhowmik and O.~Ramar\'e.
\newblock Average orders of multiplicative arithmetic functions of integer
  matrices.
\newblock {\em Acta Arithmetica}, LXVI.1:45--62, 1994.

\bibitem{bhra2}
G.~Bhowmik and O.~Ramar\'e.
\newblock On the rationality of zeta functions of finite abelian groups.
\newblock {\em Preprint}, 1997.

\bibitem{bhwu}
G.~Bhowmik and J.~Wu.
\newblock On the asymptotic behaviour of the number of subroups of finite
  abelian groups.
\newblock {\em Arch.Math}, 69:95--104, 1997.

\bibitem{BFOW}
L.~Brekke, P.G.O. Freund, M.~Olson, and E.~Witten.
\newblock Non-archimedean string dynamics.
\newblock {\em Nuclear Physics}, B302:365--402, 1988.

\bibitem{dahlquist}
G.~Dahlquist.
\newblock On the analytic continuation of eulerian products.
\newblock {\em Ark. Mat.}, 1:533--554, 1952.

\bibitem{bretecheasterisque}
R.~de~la {B}ret\`eche.
\newblock Sur le nombre de points de hauteur born\'ee d'une certaine surface
  cubique singuli\`ere.
\newblock {\em Peyre, Emmanuel (ed.), Nombre et r\'epartition de points de
  hauteur born\'ee. Ast\'erisque}, 251(2):51--57, 1998.

\bibitem{bretechetorique}
R.~de~la {B}ret\`eche.
\newblock Compter des points d'une vari\'et\'e torique.
\newblock {\em J. Number Theory}, 87(2):315--331, 2001.

\bibitem{bretechecompo}
R.~de~la {B}ret\`eche.
\newblock Estimation de sommes multiples de fonctions arithm\'etiques.
\newblock {\em Compos. Math.}, 128(3):261--298, 2001.

\bibitem{dusautoy}
M.~du~Sautoy.
\newblock Zeta functions of groups and natural boundaries.
\newblock {\em Preprint}, 2003.

\bibitem{dSG00}
M.~du~Sautoy and F.~Grunewald.
\newblock Analytic properties of zeta functions and subgroup growth.
\newblock {\em Ann. of Math.}, 152(3):793--833, 2000.

\bibitem{dSG02}
M.~du~Sautoy and F.~Grunewald.
\newblock Zeta functions of groups: zeros and friendly ghosts.
\newblock {\em Amer. J. Math.}, 124(1):1--48, 2002.

\bibitem{estermann}
T.~Estermann.
\newblock On certain functions represented by dirichlet series.
\newblock {\em Proc. London Math. Soc.}, 27:435--448, 1928.

\bibitem{fouvry}
E.~Fouvry.
\newblock Sur la hauteur des points d'une certaine surface cubique
  singuli\`ere.
\newblock {\em Ast\'erisque}, 251:31--49, 1998.

\bibitem{morozcubique}
D.R. Heath-Brown and B.Z. Moroz.
\newblock The density of rational points on the cubic surface
  $x_0^3=x_1x_2x_3$.
\newblock {\em Math. Proc. Camb. Philos. Soc.}, 125(3):385--395, 1999.

\bibitem{Hoo}
C.~Hooley.
\newblock On the number of points on a complete intersection over a finite
  field.
\newblock {\em Journal of number theory}, 38:338--358, 1991.

\bibitem{igusa}
J.-I. Igusa.
\newblock Universal {\it p}-adic zeta functions and their functional equations.
\newblock {\em Amer.J.Math}, 111:671--716, 1989.

\bibitem{lichtinduke1}
B.~Lichtin.
\newblock The asymptotics of a lattice point problem associated to a finite
  number of polynomials. {I}.
\newblock {\em Duke Math. J.}, 63(1):139--192, 1991.

\bibitem{lichtinduke2}
B.~Lichtin.
\newblock The asymptotics of a lattice point problem associated to a finite
  number of polynomials. {II}.
\newblock {\em Duke Math. J.}, 77(3):699--751, 1995.

\bibitem{lichtincomp}
B.~Lichtin.
\newblock Asymptotics determined by pairs of additive polynomials.
\newblock {\em Compositio Math.}, 107(3):233--268, 1997.

\bibitem{moroz3}
B.Z. Moroz.
\newblock On the integer points of an affine toric variety (general case).
\newblock {\em Quart. J. Math. Oxford Ser. (2)}, 50(197):37--47, 1999.

\bibitem{salberger}
P.~{S}alberger.
\newblock Tamagawa measures on universal torsors and points of bounded height
  on fano varieties.
\newblock {\em Peyre, Emmanuel (ed.), Nombre et répartition de points de
  hauteur bornée. Astérisque}, 251(2):91--258, 1998.

\bibitem{tit}
E.C. Titchmarsh.
\newblock {\em The theory of the Riemann Zeta-function.}
\newblock Clarendon Press. Oxford. Revised by D.R. Heath-Brown, 2 edition,
  1986.

\bibitem{zagier}
D.~Zagier.
\newblock Values of zeta functions and their applications in ''first european
  congress of mathematics'', vol. ii, p. 497-512.
\newblock {\em Progress in Maths.}, 120 Birkhauser, 1994.

\end{thebibliography}

\end{document}